\documentclass[12pt]{amsart}
\usepackage{amscd,amssymb,graphics}
\usepackage{mathrsfs}

\newtheorem{theorem}{Theorem}[section]
\newtheorem{corollary}[theorem]{Corollary}

\theoremstyle{definition}
\newtheorem{definition}[theorem]{Definition}
\newtheorem{question}[theorem]{Open question}
\newtheorem{conjecture}[theorem]{Conjecture}
\theoremstyle{remark}
\newtheorem{remark}[theorem]{Remark}

\newtheorem{example}[theorem]{Example}

\newcommand{\abs}[1]{\lvert#1\rvert}
\def\norm#1{\left\Vert#1\right\Vert}

\def\C {{\mathbb C}}
\def\N{{\mathbb N}}
\def\Z {{\mathbb Z}}

\def\R{{\mathbb R}}

\def\tr{{\mathrm{tr}}}
\def\e{\varepsilon}

\def\Aut{{\mbox{\rm Aut}\,}}
\def\H {{\mathscr H}}
\newcommand\B{{\mathscr B}}

\newcounter{quest}
\stepcounter{quest}
\def\under#1{{$\underline{\mbox{#1}}$}}

\begin{document}

\title[Hyperlinear and sofic groups]
{Hyperlinear and sofic groups:\\ a brief guide}

\author[V.G. Pestov]{Vladimir G. Pestov}

\address{Department of Mathematics and Statistics, 
University of Ottawa, 585 King Edward Ave., Ottawa, Ontario, Canada K1N 6N5}

\email{vpest283@uottawa.ca}


%\date{24 April 2008}

% \keywords{}
% 
\begin{abstract} 
This is an introductory survey of the emerging theory of two new classes of (discrete, countable) groups, called hyperlinear and sofic groups. They can be characterized as subgroups of metric ultraproducts of families of, respectively, unitary groups $U(n)$ and symmetric groups $S_n$, $n\in\N$. Hyperlinear groups come from theory of operator algebras (Connes' Embedding Problem), while sofic groups, introduced by Gromov, are motivated by a problem of symbolic dynamics (Gottschalk's Surjunctivity Conjecture). Open questions are numerous, in particular it is still unknown if every group is hyperlinear and/or sofic. 
\end{abstract}

\maketitle

\section{Introduction}
Relatively recently, two new classes of (discrete, countable) groups have been isolated: hyperlinear groups and sofic groups. They come from different corners of mathematics (operator algebras and symbolic dynamics, respectively), and were introduced independently from each other, but are closely related nevertheless.

Hyperlinear groups have their origin in Connes' Embedding Conjecture about von Neumann factors of type $II_1$, while sofic groups, introduced by Gromov, are motivated by Gottschalk Surjunctivity Conjecture (can a shift $A^G$ contain a proper isomorphic copy of itself, where $A$ is a finite discrete space and $G$ is a group?). 

Groups from both classes can be characterized as subgroups of metric ultraproducts of families of certain metric groups (formed in the same way as ultraproducts of Banach spaces): unitary groups of finite rank lead to hyperlinear groups, symmetric groups of finite rank to sofic groups. 

We offer an introductory guide to some of the main concepts, results, and sources of the theory, following Connes, Gromov, Benjamin Weiss, Kirchberg, Ozawa, Radulescu, Elek and Szab\'o, and others, and discuss open questions which are for the time being perhaps more numerous than the results. 

The present author hopes the survey will be of interest to mathematicians of many different backgrounds.

Still, there are good reasons to publish the paper in a journal addressed to logicians. Model and set theorists have spent more time working with ultraproducts than anyone else, and in particular there is now a well-developed model theory of metric structures \cite{BYBHU}. 
% For a proof that this experience is relevant to the kind of problems discussed in this survey, look no further than the recent preprint \cite{FPS}. 
And even if groups that we consider here are abstract (no topology), they naturally appear as subgroups of certain ``infinite-dimensional'' groups, objects some of the deep recent insights in whose structure we owe, again, to logicians, see, e.g. \cite{bergman,hjorth,KR,NVT,RS}. 

Finally, if this Introduction looks more like an abstract, it is only because the rest of the paper is nothing but an extended introduction.

\section{Ultraproducts}

\subsection{Algebraic ultraproducts}
An algebraic ultraproduct of a family $(G_\alpha)_{\alpha\in A}$ of algebraic structures with regard to an ultrafilter ${\mathcal U}$ on the index set $A$ was introduced in full generality by Jerzy \L o\'s in 1955 \cite{los}. (A prehistory of the concept is discussed in \cite{BS}, Ch. 5 and 12, where it is noted that the construction is foreshadowed by the 1930's work of G\"odel and Skolem, while Hewitt in his well-known 1948 paper \cite{hewitt} was constructing non-archimedean ordered fields by means of a procedure of which an algebraic ultrapower of $\R$ is a special case.) For instance, if the $G_\alpha$ are groups, then their ultraproduct $\left(\prod_{\alpha\in A}G_{\alpha}\right)_{\mathcal U}$ is the quotient group of the cartesian product $\prod_{\alpha\in A}G_{\alpha}$ by the normal subgroup $N_{{\mathcal U}}$ consisting of all threads $g=(g_{\alpha})$ with the property $\{\alpha\colon g_{\alpha}=e_\alpha\}\in{\mathcal U}$. (Which can be expressed in an eye-catching way by the formula $\lim_{\alpha\to{\mathcal U}}g_{\alpha}=e$.)

\subsection{Ultraproducts of normed spaces}
The above concept can be refined to suit some situations where the algebraic structures $G_{\alpha}$ possess a metric.

Historically the first such case was the ultraproduct of a family of normed spaces -- or, in the language of non-standard analysis, the nonstandard hull of an internal normed space. A particular case of a Banach space ultrapower of a single normed space (or, which is more or less the same, the nonstandard hull of a standard normed space) can be found in Abraham Robinson's {\em Nonstandard Analysis} \cite{robinson} (at the end of subsection 7.1). A general case was treated by W.A.J. Luxemburg \cite{luxemburg} (in the framework of nonstandard analysis) and, independently, by Dacunha-Castelle and Krivine \cite{DCK}. 
For a modern overview of this line of research, see the recent survey \cite{HI}, while an even more general setting of metric spaces is dealt with in the book \cite{BYBHU}.

Dusa McDuff \cite{mcduff} and, independently, Janssen \cite{janssen} had introduced ultraproducts of finite von Neumann algebras at about the same time; we will consider this construction in Section \ref{s:cec}. 

Let $E_\alpha$, $\alpha\in A$ be a family of normed spaces and let ${\mathcal U}$ be an ultrafilter on the index set $A$. Just like in the discrete case, the ultraproduct of the above family will be a quotient space, but (i) of a generally proper subspace of $\prod_{\alpha}E_\alpha$, and (ii) by a larger subspace than $N_{{\mathcal U}}$. Namely, define a normed linear space
\[{\mathscr E}=\oplus^{\ell^\infty}E_\alpha=\left\{x\in\prod_{\alpha}E_\alpha\colon \sup_{\alpha}\norm{x_\alpha}<\infty\right\},\] 
which is, in a certain sense, the largest linear subspace of the cartesian  product of $E_\alpha$'s over which one can define a norm extending the norms on $E_\alpha$:
\[\norm x = \sup_{\alpha\in A}\norm{x_\alpha}_{\alpha}.\]

Now define the subspace of ``infinitesimals,''
\[{\mathscr N}=\left\{x\colon \lim_{\alpha\to{\mathcal U}}\norm{x_{\alpha}}=0\right\}.\]
The limit along the ultrafilter is defined as the number $a$ with the property that for every $\e>0$,
\[\left\{\alpha\colon \abs{x_{\alpha}-a}<\e\right\}\in {\mathcal U}.\]
A convenient feature of this concept is that every bounded sequence of reals has an ultralimit along a given ultrafilter, which is of course a restatement of the Heine--Borel theorem, with the same proof. 
The linear subspace ${\mathscr N}={\mathscr N}_{{\mathcal U}}$ is closed in ${\mathscr E}$, and the normed space 
\[\left(\prod E_\alpha\right)_{\mathcal U} ={\mathscr E}/{\mathscr N}_{{\mathcal U}}\]
is called the Banach space ultraproduct of the family $(E_{\alpha})$ modulo the ultrafilter ${\mathcal U}$. 

Here is a ``direct'' definition of a norm on the ultraproduct:
\[\norm{\bar x} =\lim_{\alpha\to{\mathcal U}}\norm{x_{\alpha}}.\]

In the language of nonstandard analysis, the same object will be obtained by choosing an ``infinitely large'' (external) index $\nu\in\,^\ast\! A\setminus A$ and forming the quotient of the (external) linear space ${\mathrm{fin}\,}E_\nu$ of all elements with finite norm by the monad of zero, $\mu(0)$, consisting of all infinitesimals of $E_{\nu}$. The norm of a coset containing $x$ is set equal to the standard part of $\norm x$. The space obtained this way is known as the {\em nonstandard hull} of $E_{\nu}$ and denoted $\widehat{E_{\nu}}$. The freedom in choosing an external index $\nu$ corresponds to the freedom of choosing an ultrafilter $\mathcal U$ in the ultraproduct construction.

A sufficiently accurate rendering of Cantor's diagonal argument shows that if the ultrafilter ${\mathcal U}$ is not countably complete (in particular, is non-principal), then the ultraproduct $\left(\prod E_\alpha\right)_{\mathcal U}$ is a Banach space. If for some natural number $n$  the set of indices $\alpha$ with $\dim E_{\alpha}\leq n$ is in ${\mathcal U}$, the ultraproduct is of dimension $n$ itself. 
Otherwise, under the same assumption on ${\mathcal U}$, the ultraproduct is a non-separable Banach space, which is again shown through a variation of the diagonal argument. 

\subsection{Ultraproducts of metric groups: first attempt} We want 
to generalize the above construction to the case of metric groups. Let us recall that a metric $d$ on a group $G$ is {\em left-invariant} if $d(gx,gy)=d(x,y)$, for all $g,x,y\in G$. If a topological group $G$ is metrizable, then there exists a compatible left-invariant metric by the classical Kakutani theorem.
So, let $(G_\alpha,d_\alpha)$ be a family of topological groups equipped with compatible left-invariant metrics, and let ${\mathcal U}$ be an ultrafilter on $A$.
We will just emulate, word for word, the construction in the case of Banach spaces, and form, first, the ``finite part'' of the cartesian product:
\[{\mathscr G}=\left\{x\in\prod_{\alpha}G_\alpha\colon \sup_{\alpha}d(x_\alpha,e)<\infty\right\}.\]
This ${\mathscr G}$ is indeed a subgroup of the product, as follows from a simple estimate (which uses left invariance in an essential way):
\begin{eqnarray*}
d(xy,e)&=&d(y,x^{-1})\\ &\leq& d(y,e)+d(x^{-1},e) \\
&=&d(y,e)+d(e,x).
\end{eqnarray*}
The same inequality shows that
\[{\mathscr N}=\left\{x\colon \lim_{\alpha\to{\mathcal U}}d(x_{\alpha},e)=0\right\}\]
is a subgroup of ${\mathscr G}$.
However, $\mathscr N$ is not necessarily normal in $\mathscr G$: in general,
\[\lim_{\alpha\to{\mathcal U}}d(g_\alpha^{-1} x_{\alpha}g_\alpha,e)\neq 0.\]
 
Here is an example that is sufficiently interesting in itself to merit a discussion (especially given the amount of attention the infinite symmetric group has been getting from logicians recently, cf. e.g. \cite{bergman,KR}).

\begin{example}
Let $S_\infty$ denote the infinite symmetric group consisting of all self-bijections of a countably infinite set $\omega$, with its standard Polish topology induced from the embedding into the Tychonoff power $S_\infty\hookrightarrow \omega^{\omega}$, where $\omega$ is viewed as a discrete topological space. This topology is generated by 
the following left-invariant metric:
\[d(\sigma,\tau) =\sum_{i=1}^\infty \{2^{-i}\colon \sigma(i)\neq\tau(i)\}.\]
This metric can be interpreted, assuming a viewpoint of ergodic theory, as the so-called {\em uniform metric} on the group $S_\infty$ considered as the group of nonsingular transformations of the purely atomic probability measure space $\omega$ where every singleton $\{n\}$, $n=1,2,\ldots$, is assigned measure $2^{-n}$. The distance between two transformations is the measure of the set of points where they differ between themselves:
\[d(\sigma,\tau)=\mu\{n\colon \sigma(n)\neq\tau(n)\}.\]

\begin{figure}[ht]
\begin{center}
\scalebox{0.3}[0.3]{\includegraphics{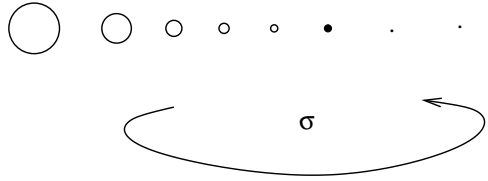}} 
\caption{$\omega$ as a purely atomic probability space.}
\label{fig:fig-dots}
\end{center}
\end{figure}

If we choose a nonprincipal ultrafilter ${\mathcal U}$ on the natural numbers, and form the subgroups $\mathscr G$ and $\mathscr N$ of $(S_\infty)^{\N}$ as above, then $\neg({\mathscr N}\triangleleft{\mathscr G})$. Indeed, ${\mathscr G}= (S_\infty)^{\N}$. Now
consider two sequences of transpositions of $\omega$, $x=(x_i)=((i,i+1))_{i\in\omega}$ and $y=(y_i)=((1,i))$. Then $d(x_i,e)=2^{-i}+2^{-(i+1)}\to 0$, so $x\in {\mathscr N}$.
% , while $d(y_i,e)<\sum_i2^{-i}=1$, so $y\in {\mathscr{G}}$.
(To verify the distance estimates, use Figure \ref{fig:fig-dots}.) At the same time, since $y_i^{-1}x_iy_i = (1,i+1)$, 
\[d(y_i^{-1}x_iy_i,e)\to 1/2.\]
\end{example}

\subsection{Bi-invariant metrics}
The message of the previous example is that in order to form ultraproducts of metric groups, it is necessary to use {\em bi-invariant} metrics:
\[d(gx,gy)=d(x,y)=d(xg,yg).\]
Such metrics on groups always determine the so-called SIN topologies, i.e. topologies for which left and right uniformities coincide, or --- equivalently --- open subsets invariant under conjugation form a neighbourhood basis of identity. (Hence the acronym: Small Invariant Neighbourhoods.)

If $(G_\alpha,d_{\alpha})$ is a family of groups equipped with bi-invariant metrics and ${\mathcal U}$ is an ultrafilter on the index set $A$, then the subgroup $\mathscr N$ of ``infinitesimals'' is normal in the subgroup $\mathscr G$ of finite elements, and the quotient group 
\[\left(\prod_{\alpha\in A}G_{\alpha}\right)_{\mathcal U}
={\mathscr G}/{\mathscr N}\]
is well-defined. Equipped with the bi-invariant metric 
\[d(x{\mathscr N},y{\mathscr N})=\lim_{\alpha\to{\mathcal U}}d_{\alpha}(x_{\alpha},y_{\alpha})\]
(and the corresponding group topology), it will be referred to as the {\em metric ultraproduct} of the family $(G_\alpha,d_{\alpha})$ modulo ${\mathcal U}$. 

Just as in the case of normed spaces, the ultraproduct of a family of groups with bi-invariant metrics is a complete topological group, which is either non-separable or locally compact (of course assuming ${\mathcal U}$ to be non countably complete).

Moreover, in all the examples we will be considering below, $\mathscr G$ coincides with the full cartesian product, because all the metrics are uniformly bounded from above. (In fact, one can always replace a bi-invariant metric $d$ on a group with e.g. the bi-invariant metric $\min\{d,1\}$, so this is not much of an issue.) 

Here are a few of the most important examples of groups equipped with natural bi-invariant metrics.

\begin{example}
The {\em uniform metric} on $\Aut(X,\mu)$, the group of measure-preserving transformations of a finite measure space $(X,\mu)$:
\[d(\sigma,\tau)=\mu\{x\in X\colon \sigma(x)\neq \tau(x)\}.\]
\end{example}

% We have already seen that in the case of a (bigger) group of {\em measure class preserving} transformations, the uniform metric need not be invariant on the right. However, in the {\em measure preserving} case bi-invariance follows easily.

A particular case of the above construction is:

\begin{example}
The {\em normalized Hamming distance} on the symmetric group $S_n$ of finite rank $n$ is given by
\[d_{hamm}(\sigma,\tau)=\frac 1n\sharp\left\{i\colon\sigma(i)\neq\tau(i)\right\}.\]
The measure space in question is a finite set $[n]=\{1,2,\ldots,n\}$, equipped with the uniform ($=$ normalized counting) measure: $\mu\{i\}=1/n$ for every $i$.
\end{example}
 
\begin{example}
Let $\H$ be a Hilbert space (either finite or infinite dimensional). Denote by $U(\H)$ the group of all unitary operators on $\H$,
\[U(\H)=\{u\colon\H\to\H\mid \mbox{$u$ is linear and bounded and }u^\ast u=uu^\ast ={\mathrm{Id}}\},\]
and equip it with the {\em uniform operator metric:} 
\[d_{unif}(u,v)=\norm{u-v}=\sup_{\norm x\leq 1}\norm{(u-v)(x)}.\]
This metric is easily checked to be bi-invariant (and the topology it induces is known as the {\em uniform operator topology}). 
\end{example}

\begin{remark}
\label{r:kirchberg}
Sometimes properties of ultraproducts depend on non-principal ultrafilters with regard to which the ultraproducts are formed. 
Though the following question is only indirectly linked to the topic of these notes (with papers \cite{mcduff} and \cite{connes-injective} providing a link), it illustrates the point.

Notice that every (metric) group $G$ embeds into its own ultrapower $\left(G^I\right)_{\mathcal U}$ diagonally under the map $x\mapsto (x,x,\ldots)$ as a metric subgroup, and recall that for every Hilbert space $\H$ the centre of $U(\H)$ is the circle group $\{\lambda\cdot{\mathrm{Id}}\colon \abs\lambda=1\}$.

\begin{question}[Kirchberg, cf. \cite{kirchberg06}, question 2.22 on p. 195] Is the centralizer of the subgroup $U(\ell^2)$ in the metric ultrapower
\[\left((U(\ell^2),d_{unif})^{\N}\right)_{\mathcal U}\]
equal to $\{\lambda\cdot{\mathrm{Id}}\colon \abs\lambda=1\}$?
\end{question}

Ilijas Farah and N. Christopher Phillips have shown that there is always an ultrafilter $\mathcal U$ for which the centralizer is nontrivial (cf. \cite{farah1}, an article is in preparation). It remains unknown if ultrafilters for which the centralizer is trivial do exist. 
% They conjecture that possibly the existence of an ultrafilter $\mathcal U$ for which the centralizer is trivial is independent from ZFC.
\end{remark}

\begin{example}
If $\H=\C^n$ is $n$-dimensional Hermitian space, then the group $U(\H)$ is denoted $U(n)$ and called the unitary group of rank $n$. It can be identified with the group of all $n\times n$ unitary matrices with complex entries, $u=(u_{ij})_{i,j=1}^n$.
%The following will be of particular importance to us. 

The {\em normalized Hilbert-Schmidt metric} on the group $U(n)$ is the standard $\ell^2$ distance between matrices viewed as elements of an $n^2$-dimensional Hermitian space, which is normalized so as to make the identity matrix have norm one:
\[d_{HS}(u,v)=\norm{u-v}_2=\sqrt{\frac 1n\sum_{i,j=1}^n \abs{u_{ij}-v_{ij}}^2}.
\]
In order to verify bi-invariance of the Hilbert-Schmidt distance, it suffices to rewrite it using trace of a matrix, as follows:
\[d_{HS}(u,v)= \frac 1{\sqrt n}\sqrt{{\mathrm{tr}\,} ((u_{ij}-v_{ij})^\ast(u_{ij}-v_{ij}))}.\]
In this form, bi-invariance follows from the main property of trace: ${\mathrm{tr}}\,(AB)={\mathrm{tr}}\,(BA)$.
\end{example}

\section{Definitions}

Now hyperlinear and sofic groups can be defined side by side, in a completely analogous fashion.

\begin{definition}
A group $G$ is {\em sofic} if it is isomorphic to a subgroup of a metric ultraproduct of a suitable family of symmetric groups of finite rank with their normalized Hamming distances. In other words, there are a set $A$, an ultrafilter ${\mathcal U}$ on $A$, and a mapping $\alpha\mapsto n({\alpha})$ so that 
\begin{equation}
\label{eq:sofic}
G<\left(\prod_{\alpha} (S_{n({\alpha})},d_{hamm})\right)_{\mathcal U}.\end{equation}
\end{definition}

\begin{definition}
A group $G$ is {\em hyperlinear} if it is isomorphic to a subgroup of a metric ultraproduct of a suitable family of unitary groups of finite rank, with their normalized Hilbert-Schmidt distances. In other words, there are a set $A$, an ultrafilter ${\mathcal U}$ on $A$, and a mapping $\alpha\mapsto n({\alpha})$ so that 
\begin{equation}
G<\left(\prod_{\alpha} (U(n(\alpha)),d_{HS})\right)_{\mathcal U}.\end{equation}
\end{definition}

Perhaps the most natural first question that comes to mind, is:
what is the relation between the two classes of groups which seem to be so similar?
The finite permutation group $S_n$ embeds into the unitary group $U(n)$ as a subgroup if we associate to a permutation $\sigma$ the corresponding permutation matrix $A_{\sigma}$ the way we do it in a second year Linear Algebra course:
\[\left(A_{\sigma}\right)_{ij}=\begin{cases} 1,&\mbox{ if }\sigma(j)=i,
\\
0,&\mbox{ otherwise.}
\end{cases}
\]

One has to be careful here: the restriction of the normalized Hilbert-Schmidt distance to $S_n$ is not, in fact, even Lipschitz equivalent to the normalized Hamming distance. Nevertheless, the two distances agree with each other sufficiently well so as to preserve the embeddings at the level of metric ultraproducts and lead to the following result.  

\begin{theorem}[Elek and Szab\'o \cite{ES}] Every sofic group is hyperlinear. 
\end{theorem}

\begin{proof}
Let us compare the values of two distances (the normalized Hamming distance and the normalized Hilbert-Schmidt distance) between two permutations, $\sigma,\tau\in S_n$:
\begin{eqnarray*}
d_{hamm}(\sigma,\tau) &=& \frac 1n\sharp\{i\colon\sigma(i)\neq\tau(i)\}\\
&=&  \frac 1n\sharp\{i\colon(\sigma\tau^{-1})(i)\neq i\}\\
&=& \frac 1 {2n} {\mathrm{tr}}\,({\mathbb{I}}-A_{\sigma\tau^{-1}})+\frac 1 {2n} {\mathrm{tr}}\,({\mathbb{I}}-A_{\tau\sigma^{-1}}) \\
&=& \frac 1 {2n} {\mathrm{tr}}\,((A_{\sigma}-A_\tau)(A_\sigma-A_\tau)^\ast) \\
&=&\frac 12 \left(d_{HS}(A_\sigma,A_\tau)\right)^2.
\end{eqnarray*}
We conclude: the condition 
\[d(x_n,e)\to_{\mathcal U} 0\] 
is the same with regard to both metrics, so as topological groups, the metric ultraproduct of $S_n$'s embeds into the metric ultraproduct of $U(n)$'s over the same ultrafilter:
\begin{eqnarray*}
\left(\prod_{\alpha} S_{n(\alpha)}\right)_{\mathcal U}&=&\prod S_{n(\alpha)}/ \left({\mathscr N}\cap \prod S_{n(\alpha)}\right)\\[4mm] &<& \prod_{\alpha} U(n(\alpha))/ {\mathscr N}\\
&=&\left(\prod_{\alpha} U(n(\alpha))\right)_{\mathcal U}.
\end{eqnarray*}
\end{proof}

By contrast, the converse implication is unknown.

\begin{question}
Is every hyperlinear group sofic? 
\end{question}

Here one can speculate that since, by Malcev's theorem (\cite{malcev}, also cf. Theorem 6.4.13 in \cite{BO}), every finitely generated subgroup of $U(n)$ is residually finite and thence sofic (Example \ref{ex:resfin} below), a likely answer might be ``yes,'' but this remains just this author's guess.

Bearing in mind Remark \ref{r:kirchberg}, we will address the following question: to what extent do the two concepts depend on the choice of an ultrafilter? To this end, we will reformulate both definitions in a way not using ultraproducts.
Here is an equivalent reformulation of the concept of a sofic group.

\begin{theorem}
\label{th:soficcriterion}
A group $G$ is sofic if and only if for every finite $F\subseteq G$ and for each  $\e>0$, there exist a natural $n$ and a mapping $\theta\colon F\to S_n$  so that 
\begin{enumerate}
\item if $g,h,gh\in F$, then $d_{hamm}(\theta(g)\theta(h), \theta(gh))<\e$,
\item if $e\in F$ then $d(\theta(e),{\mathrm{Id}})<\e$, and 
\item for all distinct $x,y\in F$, $d_{hamm}(\theta(x),\theta(y))\geq 1/4$. 
\end{enumerate}
{\rm (}A mapping satisfying conditions (1)-(2) is called an $(F,\e)${\em -almost homomorphism}{\rm .)}
\end{theorem}

\begin{proof}
$\Leftarrow$: A family of $(F,1/n)$-homomorphisms $\theta_{F,1/n}$ from a group $G$ to symmetric groups $S_{(F,1/n)}$ of finite rank satisfying condition (3)  determines a group monomorphism $G\hookrightarrow \left(\prod S_{(F,1/n)}\right)_{\mathcal U}$ in the following standard manner. First choose as the index set $A$ the collection of all pairs $(F,1/n)$, where $F$ is a finite subset of $G$ and $n\geq 1$, partially ordered in a natural way. Next choose a nonprincipal ultrafilter $\mathcal U$ on $A$ which contains every subset of the form $\left\{(\Phi,1/m)\colon \Phi\supseteq F,m\geq n\right\}$. Then the mapping 
\begin{eqnarray*}
\theta\colon G&\to& \left(\prod_{(F,1/n)}S_{(F,1/n)}\right)_{\mathcal U}, \\
g&\mapsto & \left(\theta_{F,1/n}(g)\right)_{\mathcal U}
\end{eqnarray*}
is a well-defined group homomorphism which satisfies $d(\theta(g),\theta(h))\geq 1/4$ for all $g\neq h$ and hence is a monomorphism.
\\[3mm]
$\Rightarrow$: If $G$ is a sofic group and $\theta$ is an embedding of $G$ into the ultraproduct as in Eq. (\ref{eq:sofic}), then a family $\theta_{\alpha}\colon G\to S_{n(\alpha)}$ of $(F,\e)$-almost homomorphisms is obtained in a routine way by considering the coordinate projections. The third condition is a little bit less straightforward, because all one can claim, is that for $g\neq h$, the images of $g$ and $h$ in the ultraproduct are distinct, but they certainly do not need to be at a distance $\geq 1/4$ or anything of the kind. To achieve the desired separation between the images of two given elements, one employs a trick known in some areas of functional analysis as ``amplification.'' 

\begin{figure}[ht]
\begin{center}
\scalebox{0.35}[0.35]{\includegraphics{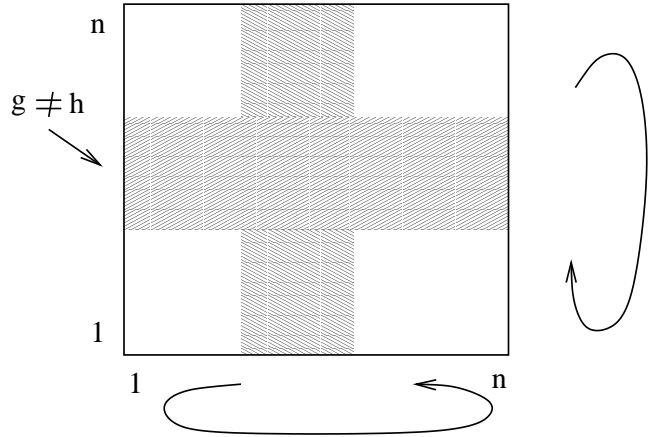}} 
\caption{The amplification trick for $k=2$.}
\label{fig:fig-tensor}
\end{center}
\end{figure}

Let $g,h\in F$, $g\neq h$. Fix an $(F,\e^\prime)$-almost homomorphism $\theta\colon F\to S_n$ with $\theta(g)\neq\theta(h)$, where $\e^\prime>0$ is sufficiently small, to be specified later. 
Denote $\delta=d(\theta(g),\theta(h))=\mu_{\sharp}\{i\colon \theta(g)(i)\neq\theta(h)(i)\}$, where $\mu_{\sharp}$ denotes the normalized counting measure on $[n]$. All we can claim, is that $\delta>0$. 
Now ``amplify'' $\theta$ by re-embedding $F$ diagonally into the group of permutations of the square $[n]\times [n]=[n^2]$:
\[\theta^{\otimes 2}(g) = (\theta(g),\theta(g)),\]
where the latter acts on $[n]\times [n]$ by double permutations. The measure of the set of pairs $(i,j)$ on which $(\theta(g),\theta(g))$ and $(\theta(h),\theta(h))$ are different (a cross in Fig. \ref{fig:fig-tensor}) has increased from $\delta$ to at least $1-(1-\delta)^2=2\delta-\delta^2$. Amplifying $\theta$ as many times as necessary, one can increase the distance between any given pair of points to $1/4$. In fact, the choice of $1/4$ here was arbitrary, and one can replace $1/4$ with any real number strictly between $0$ and $1$.

Finally, let us address the choice of $\e^\prime>0$. A somewhat undesirable outcome of amplification is that the distances between $\theta(x)\theta(y)$ and $\theta(xy)$ will also increase so that $\theta^{\otimes 2}$ is no longer an $\e^\prime$-almost homomorphism. However, the remedy here is simple. Since the value of $\delta$ can be assumed as close to $d(g,h)$ as we wish, the desired number $k$ of amplifications can be estimated before $\theta$ is chosen. 
One starts with an $\e^\prime$ small enough so that $\theta^{\otimes k}$ remains an $\e$-homomorphism for a prescribed value $\e>0$. For instance, $\e^\prime=2^{-k}\e$ will do.
\end{proof}

The above result is due to Elek and Szab\'o \cite{ES}. Historically, their argument followed a similar result for hyperlinear groups, appearing in Radulescu \cite{radulescu00}:

\begin{theorem}
A group $G$ is hyperlinear if and only if for every finite 
$F\subseteq G$ and each $\e>0$ there exist a natural $n$ and a mapping $\theta\colon F\to U(n)$ (an $(F,\e)$-almost homomorphism) so that 
\begin{enumerate}
\item  if $g,h,gh\in F$, then $\norm{\theta(g)\theta(h)-\theta(gh)}_2<\e $,
\item if $e\in F$ then $\norm{\theta(e)-{\mathrm{Id}}}_2<\e$, and
\item for all distinct $x,y\in F$, $\norm{\theta(x)-\theta(y)}_2\geq 1/4$.
\end{enumerate}
\qed
\label{th:hypercriterion}
\end{theorem}

\begin{remark}
\label{r:sqrt2}
Again, the choice of $1/4$ here is totally arbitrary, and
the condition (3) can be refined so as to require $\norm{\theta(x)-\theta(y)}_2$ to be as close to $\sqrt 2$ as desired.
\end{remark}

In view of the two preceding results, 
the concepts of a hyperlinear and of a sofic group esentially do not depend on the choice of an ultrafilter. 
A countable group is hyperlinear (sofic) if and only if it embeds, as a subgroup, into the metric ultraproduct of the family $S_n$ (resp., $U(n)$), $n\in\N$, with regard to {\em some} (equivalently: {\em any}) nonprincipal ultrafilter on the natural numbers. This follows from the two previous theorems supplemented by a simple argument along the same lines as the proof of necessity ($\Rightarrow$) in Theorem \ref{th:soficcriterion}. And an apparent greater generality of allowing uncountable groups is an illusion: as follows from the two preceding results, a group $G$ is hyperlinear (sofic) if and only if so are all finitely generated subgroups of it.

Here are two central open questions of the theory.

\begin{question}
\label{q:sofic?} Is every group sofic? 
\end{question}

This question originated in 1999 Gromov's article \cite{gromov99} where the concept of a sofic group was first introduced (without a name of its own) in order to attack Gottschalk's Surjunctivity Conjecture (\cite{gottschalk}, see Conjecture \ref{co:gott} below). The current expression ``sofic group'' was coined by Benjy Weiss \cite{weiss}. 

\begin{question}
\label{q:hyperlinear?}
Is every group hyperlinear? 
\end{question}

The statement that every group (equivalently: every countable group) is hyperlinear is known as {\em Connes' Embedding Conjecture for Groups,} and we will discuss it below. The origin of this conjecture is Connes' 1979 paper \cite{connes-injective}. The expression ``hyperlinear group'' belongs to Radulescu \cite{radulescu00}.

Both questions are equivalent to their versions for countable groups, by force of Theorems \ref{th:soficcriterion} and \ref{th:hypercriterion}.

We will discuss the origins and significance of both concepts in greater detail below.

\begin{remark}
In the definition of a hyperlinear group, one can replace the unitary groups $U(n)$ by the orthogonal groups $O(n)$ (with the Hilbert-Schmidt distance), or the symplectic groups $Sp(n)$, without changing the notion.
\end{remark}

Moreover, instead of the ultraproduct of a sequence of groups, one can consider an ultrapower of a single group, by using the followng result, which follows more or less directly from Theorems \ref{th:soficcriterion} and \ref{th:hypercriterion}.

\begin{theorem}
Let $G$ be a group equipped with a bi-invariant metric $d$ and containing an increasing chain of subgroups isomorphic to $U(n)$ (respectively, to $S_n$), $n\in\N$, whose union is dense in $G$ and such that the restriction of $d$ to $U(n)$ (respectively, to $S_n$) is the normalized Hilbert-Schmidt distance (respectively, the normalized Hamming distance). Then a group $\Gamma$ is hyperlinear (respectively, sofic) if and only if $\Gamma$ embeds as a subgroup into a suitable ultrapower of $G$. 
\label{th:onegroup}
\end{theorem}

% \begin{theorem}
% Indeed, on the one hand, a metric ultrapower of $U(R)$ contains the corresponding metric ultraproduct of groups $U(n)$ because $U(n)<U(R)$. On the other hand, every group $G$ that embeds into $U(R^{\omega})$ as a subgroup admits a family of almost homomorphisms into $U(R)$. As $R$ is AFD, its unitary group contains a chain of unitary groups of increasing finite rank with dense union. Using this approximation, one can further construct a family of almost homomorphisms from $G$ to $U(n)$, $n\in\N$.
% \end{theorem}

Here are two examples of groups as in Theorem \ref{th:onegroup}.

\begin{example}
Let $\mathscr R$ be a Borel equivalence relation on a standard Borel space $X$ equipped with a finite measure $\mu$. 
The {\em full group} of $\mathscr R$ in the sense of Dye \cite{Dye}, denoted $[{\mathscr R}]$, is the subgroup of all measure class preserving transformations $\sigma$ of $(X,\mu)$ with the property $(x,\sigma(x))\in {\mathscr R}$ for $\mu$-a.e. $x$. If equipped with the uniform metric, $[{\mathscr R}]$ is a Polish group. One example of such an equivalence relation is the {\em tail equivalence relation} on $\{0,1\}^\omega$, where two infinite binary strings are equivalent if and only if they coincide at all but finitely many coordinates and the measure is the product measure. (For more, see e.g. the book by Kechris and Miller \cite{KM}.)
One can show that finite permutation groups $S_n$ can be embedded into $[{\mathscr R}]$ in such a way that their union is dense and the uniform metric induces the normalized Hamming distance on $S_n$. Thus, (countable) sofic groups are exactly countable subgroups of the metric ultrapower $\left([{\mathscr R}]^{\N}\right)_{\mathcal U}$
of the full group of the tail equivalence relation, where $\mathcal U$ is a nonprincipal ultrafilter on the natural numbers.
\end{example}

\begin{example}
The group 
\[U(\infty)_2=\{u\in U(\ell^2)\colon \norm{u-{\mathrm{I}}}_2<\infty\}\]
of Schatten class $2$ perturbations of the identity, equipped with the Hilbert-Schmidt metric is approximated in a similar way with unitary groups $U(n)$ of finite rank. Countable hyperlinear groups are exactly all countable subgroups of the ultrapower $\left(U(\infty)_2^{\N}\right)_{\mathcal U}$ formed with regard to some (any) nonprincipal ultrafilter $\mathcal U$ on the natural numbers. 

Another example of a group $G$ with the same property will appear in Rem. \ref{r:R}.
\label{ex:uinfty2}
\end{example}

\begin{remark}
It is worth pointing out that not every countable group is isomorphic to a subgroup of either $U(\infty)_2$ or $[{\mathscr R}]$, or some such group $G$ satisfying the assumptions of Theorem \ref{th:onegroup}. Namely, one can prove, using results of \cite{Bek1}, that if a group $\Gamma$ has property $(T)$ and is contained in $G$ as a subgroup, then $\Gamma$ is residually finite, cf. a similar argument in \cite{valette}.
\end{remark}

\section{Examples}

It appears that all the presently known examples of hyperlinear groups are at the same time known to be sofic.

\begin{example}
Every \under{finite} group is sofic.
\\[2mm]
$\triangleleft$
Indeed, a finite group is contained in some $S_n$ as a subgroup. \hfill
$\triangleright$
\end{example}

Recall that a group $G$ is {\em residually finite} if it admits a separating family of homomorphisms into finite groups.

\begin{example}
\label{ex:resfin}
Every \under{residually finite} group is sofic.
\\[2mm]
$\triangleleft$
If $F\subseteq G$ is finite, there exists a normal $N\triangleleft G$ with $N\cap (F\cdot F^{-1})\subseteq \{e\}$ and $G/N$ finite.
The composition of the quotient mapping
$h\colon G\to G/N$ with a standard embedding of $G/N$ into the symmetric group $S_{\abs{G/N}}$ of the set $G/N$
is a $(F,\e)$-homomorphism (for all $\e>0$), satisfying the condition (3) from Theorem \ref{th:soficcriterion}.\hfill $\triangleright$
\end{example}

In particular,

\begin{example}
Every nonabelian \under{free} group (e.g. $F_\infty$) is sofic.
\end{example}

(There are different proofs of residual finiteness of free groups. A  beautiful argument of Sanov \cite{sanov}, which the present author has learned from \cite{dlH3}, pp. 25--26,
embeds $F_2$ as a subgroup into the --- obviously residually finite --- group ${\mathrm{SL}}(2,\Z)$. For another proof and historical references, as well as more examples of residually finite groups, see \cite{MKS}, pp. 116 and 414.)

Hyperlinearity of nonabelian free groups, established in 1976 independently by Connes \cite{connes-injective} and S. Wassermann \cite{wassermann}, marked the beginning of the present direction of research. In all the fairness, the proof of the result was not quite so easy as might be suggested by mere Example \ref{ex:resfin}, because it also in essence included what has later become Theorems \ref{th:hypercriterion} and \ref{th:iff}.
 
\begin{example}
Every \under{amenable} group is sofic. 
\\[2mm]
$\triangleleft$
It is plausible that of all mathematical concepts, amenability of a group admits the largest known number of equivalent definitions. (``Approriximately $10^{10^{10}}$,'' according to the authors of \cite{BO}, p. 48, where also a brief introduction to the concept can be found. For more detailed references, see \cite{Gre1,wagon}.)
One of the best known among those equivalent definitions, the {\em F\o lner condition,} says the following.
Given a finite $F\subseteq G$ and $\e>0$, there is a finite $\Phi\subseteq G$ (a {\em F\o lner set} for $F$ and $\e$) such that for each $g\in F$,
\[\vert g\Phi\bigtriangleup\Phi\vert< \e \vert\Phi\vert,\]
where $\bigtriangleup$ stands for the symmetric difference. 
The map $x\mapsto gx$ is well-defined on a subset of $\Phi$ having normalizing counting measure $>1-\e$, and by extending it over the rest of $\Phi$ in an arbitrary manner so as to get a bijection, one obtains a $(F,2\e)$-homomorphism to the symmetric group $S_{\abs\Phi}$ satisfying condition (3). 
\hfill
$\triangleright$
\end{example}

Gromov \cite{gromov99}, p. 133, calls a group $G$ {\em initially subamenable} if, given a finite subset $F\subseteq G$, one can find an amenable group containing a copy of $F$, with the same partial multiplication. In other words, one cannot tell $G$ apart from an amenable group by looking at any finite piece of $G$.

For instance, every residually finite or, more generally, residually amenable group is initially subamenable. (A group is residually amenable if homomorphisms to amenable groups separate points.)
So is every {\em LEF group} $G$ in the sense of Vershik and Gordon \cite{VG}, defined by the property that one can embed every finite $F\subseteq G$ into a suitable finite group so as to preserve the partial multiplication. 

\begin{example}
Every \under{initially subamenable} group is sofic.
\\[2mm]
$\triangleleft$ 
The proof is quite clear, because soficity is a local property!
\hfill $\triangleright$
\end{example}

Here is an example of an initially subamenable group.

\begin{example}
\label{ex:baumslag}
The {\em Baumslag--Solitar group}, given by
$\langle a,b\mid ab^3a^{-1}=b^2\rangle$, is residually solvable (homomorphisms to solvable groups separate points) \cite{kropholler}, hence residually amenable and initially subamenable, in particular sofic. (This group is known to be non-residually finite.) 
\end{example}

The hyperlinearity of the Baumslag--Solitar group was first established by Radulescu \cite{radulescu00} in a difficult proof, while the above argument I learned from Goulnara Arzhantseva.

\begin{remark}
An example of an initially subamenable (even LEF) group that is not residually amenable has been constructed in \cite{ES2}, refining a construction from \cite{VG}.
\end{remark}

The following was pointed out to me, independently, by Denis Osin and by Simon Thomas.

\begin{remark}
\label{r:osin-thomas}
Not every group is initially subamenable. 

For example, every finitely
presented non-amenable simple group $G$ (cf. \cite{BM}) is not in this class. Indeed, if a group $H$ is generated by a set having the same partial multiplication table as a sufficiently large ball in $G$, then $H$ is isomorphic to $G$ since $G$ is finitely presented and simple. 
\end{remark}

\begin{question}
May it happen that a group $G$ is sofic without being initially subamenable?
\end{question}

In Gromov's opinion (\cite{gromov99}, p. 157, line ${-14}$), ``it may (?) happen''. As noticed by Simon Thomas, if any of the known examples of finitely presented simple non-amenable groups is verified to be sofic, it will provide a distinguishing example in view of Remark \ref{r:osin-thomas}. Existence of countable simple {\em finitely generated} non-amenable sofic groups was established in \cite{ES} (cf. Corollary 3.1).

Example \ref{ex:baumslag} suggests:

\begin{question}
Is every one-relator group (that is, a group defined by a single defining relation between generators) hyperlinear (Nate Brown)? Sofic?
\label{q:1rel?}
\end{question}

Notice that

\begin{question}
Is every finitely presented group  hyperlinear (sofic)?
\end{question}
\noindent
is already equivalent to the question about arbitrary groups. Indeed, if $F$ is a finite subset of a group $G$ represented as a quotient of a free group $F_n$ by a normal subgroup $N\triangleleft F_n$, there is a finitely generated normal subgroup $N^\prime\triangleleft F_n$, $N^\prime\subseteq N$, so that $F$ embeds into $F_n/N^\prime$ with the same partial multiplication.

\section{Further criteria of soficity}

A (directed) graph $\Gamma$ is {\em edge-coloured} if there is a set $C$ (whose elements are viewed as {\em colours}) and a mapping associating to every (directed) edge $(x,y)\in E(\Gamma)$ an element of $C$. In this case, we will also say that $\Gamma$ is {\em edge $C$-coloured}.

Let now $G$ be a finitely-generated group. Fix a finite symmetric set $V$ of generators of $G$ not containing the identity $e$. The {\em Cayley graph} of $G$ (corresponding to $V$) has $G$ as set of vertices, with $(g,h)$ being adjacent if and only if $g^{-1}h\in V$, that is, if one can get to $h$ by multiplying $g$ with a generator $v\in V$ on the right. Thus, it is a {\em directed} graph.
Clearly, such a $v$ associated to an edge is unique, and so the Cayley graph of $G$ becomes an edge $V$-coloured graph.

Speaking of an $N$-ball in a (connected) graph $\Gamma$, we will mean a closed ball of radius $N$ with regard to the path distance, i.e., the length of the shortest path between two vertices. In case of a group, the path distance becomes the left-invariant {\em word distance}, $d_V(g,h)=d_V(g^{-1}h,e)$, that is, the length of the shortest $V$-word representing $g^{-1}h$. 

\begin{theorem}[Elek and Szab\'o \cite{ES}]
Let $G$ be a group with a finite generating set $V$. 
Then $G$ is sofic if and only if the following Gromov's condition is satisfied:
\\[2mm]
{\rm ({\large $\star$})}
for every natural $N$ and $\e>0$ there is a finite edge $V$-coloured graph $\Gamma$ with the property that for the fraction of at least $(1-\e)\abs \Gamma$ of vertices $x$ of $\Gamma$ the $N$-ball $B_N$ around $x$ is isomorphic, as an edge $V$-coloured graph, to the $N$-ball in $G$. 
\label{th:graph}
\end{theorem}

The condition {\rm ({\large $\star$})} says that locally $\Gamma$ looks like $G$ at all but $<\e\vert\Gamma\vert$ of its vertices. In other words, one can cut out of the Cayley graph of $G$ sufficiently many copies of the ball $B_N(e)$ and glue them together in such a fashion that most vertices of the resulting edge-coloured graph $\Gamma$ are centres of the pasted balls.

Sofic groups were originally introduced by Gromov in \cite{gromov99}, on p. 157 namely in the form of the condition ({\large $\star$}).

\begin{proof}[Proof of Theorem \ref{th:graph}, sketch]
$\Rightarrow$: The graphs $\Gamma$ are obtained by tinkering with Cayley graphs of those permutation groups $S_n$ where $G$ is mapped almost homomorphically.
\\[2mm]
$\Leftarrow$: in the presence of an edge-colouring, every element $w\in B_{N}$ determines a unique translation of $\Gamma$ that is well-defined at all but $<\e\abs\Gamma$ of its vertices (just follow any particular sequence of colours leading up to $w$ in the original ball). This way, almost homomorphisms of $G$ into finite permutation groups are constructed.
\end{proof}

Recall that an action of a group $G$ on a set $X$ is {\em free} if for every $g\in G$, $g\neq e$, and each $x\in X$ one has $gx\neq x$. If $X$ is a space with (finitely-additive) measure $\mu$, defined on some algebra of sets, both the notion of an action and that of freeness can be weakened. A {\em near-action} of $G$ on $(X,\mu)$ is an assignment to every $g\in G$ of a measure-preserving map $\theta(g)\colon X\to X$ defined $\mu$-almost everywhere, in such a way that for every $g,h\in G$, one has $\theta(g^{-1})\theta(h)=\theta(g^{-1}h)$ $\mu$-a.e. A near-action is {\em essentially free} if for every $g\in G\setminus\{e\}$, $\theta(g)(x)\neq x$ $\mu$-a.e. 

\begin{theorem}[Elek and Szab\'o \cite{ES}]
A group $G$ is sofic if and only if it admits an essentially free near-action on a set $X$ equipped with a finitely-additive probability measure $\mu$ defined on the family ${\mathscr P}(X)$ of all subsets of $X$.
\label{th:essfree}
\end{theorem}

\begin{proof}[Sketch of the proof]
$\Rightarrow$: Let $G$ be a sofic group. Choose a complete family of $(F,1/k)$-almost homomorphisms $\theta_{F,k}\colon G\to S_{n(F,k)}$ satisfying condition (3) from Theorem \ref{th:soficcriterion}, where however for every $k$ the constant $1/4$ is replaced with a constant $c_k<1$ in such a way that $c_k\to 1$ as $k\to\infty$. Each $\theta_{F,k}$ defines an ``almost action'' of $F$ on the finite set $[n]=[n(F,k)]$. Define a set
\[X=\cup_{F,k}[n(F,k)].\]
Let ${\mathcal U}$ be any ultrafilter on the directed set of all pairs $(F,k)$ containing every upper cone of the form $\{(F,k)\colon F\supseteq F_0,~k\geq k_0\}$. Now define a finitely-additive probability measure $\mu$ on the power set of $X$ by letting for each $A\subseteq X$
\[\mu(A)=\lim_{(F,k)\to{\mathcal U}}\frac{\abs{A\cap [n(F,k)]}}{n(F,k)}.\]
It is easy to see that every $g\in G$ defines $\mu$-a.e. a measure-preserving transformation of $X$, and that the near-acion of $G$ on $(X,\mu)$ defined in this way satisfies all the required properties.
\\[3mm]
$\Leftarrow$: Here one resorts to a suitable modification of the technique of paradoxical decompositions in order to exclude the existence of an invariant finitely additive measure for near-actions of non-sofic groups.
\end{proof}

The paradoxical decompositions of Banach and Tarski had motivated the very concept of an amenable group, cf. e.g. \cite{wagon}, which stresses yet again that soficity is a younger sister of amenability.

Open problem \ref{q:sofic?} is therefore equivalent to:

\begin{question}
Does every countable group admit a near action as above?
\end{question}

Theorem \ref{th:essfree} should be compared to a known characterization of amenable discrete countable groups $G$ as those admitting a left-invariant finitely-additive measure defined on all subsets of $G$. 
% Thus, yet again: soficity is a weakening of amenability in a rather natural way.

It is known, for example, that every countable discrete group $G$ acts freely on a Cantor space admitting an invariant sigma-additive Borel probability measure, $\mu$, moreover there is an explicit construction of such an action \cite{HM}. However, this $\mu$ is only defined on Borel subsets of $X$.

\section{Gottschalk Surjunctivity Conjecture}

Let $G$ be a (countable) group, $A$ a finite set equipped with discrete topology. The Tychonoff power $A^G$ is a Cantor space (i.e., a compact metrizable zero-dimensional space without isolated points), upon which $G$ acts by translations:
\[(g\cdot x)_h=x_{g^{-1}h}.\]
Equipped with this action of $G$ by homeomorphisms, $A^G$ is a {\em symbolic dynamical system,} or a {\em shift.}

\begin{conjecture}[Gottschalk Surjunctivity Conjecture, \cite{gottschalk}, 1973]
For every countable group $G$ and every finite set $A$, the shift system $A^G$ contains no proper closed $G$-invariant subspace $X$ isomorphic to $A^G$ itself (as a compact $G$-space).
\label{co:gott}
\end{conjecture}

It seems to be unclear whether it suffices to set $A=\{0,1\}$.

\begin{question}
Is the Gottschalk Surjunctivity Conjecture equivalent to its particular case where $A=\{0,1\}$?
\end{question}

Here is the main advance to date.

\begin{theorem}[Gromov \cite{gromov99}]
\label{th:gromov}
The Gottschalk Surjunctivity Conjecture holds for sofic groups.
In other words, if $G$ is a sofic group and $A$ a finite set, then the compact $G$-space $A^G$ contains no proper isomorphic copies of itself.
\end{theorem}

Now the significance of open problem \ref{q:sofic?} becomes clear: if every group is sofic, then both the Gottschalk Surjunctivity Conjecture and the Connes Embedding Conjecture for groups are true. Conversely, if the Surjunctivity Conjecture is disproved, it would imply the existence of groups that are non-sofic (though perhaps not necessarily non-hyperlinear).
\vskip 1.5mm

We will now look at some special cases of Theorem \ref{th:gromov}.

\subsubsection*{Case 1: $G$ is finite}$\,$ \\[1mm]
Here the proof is obvious for cardinality reasons, as $A^G$ is itself finite.

\subsubsection*{Case 2: $G=\Z$}$\,$ \\[1mm]
For each $n$, the set of $n$-periodic points is finite, and every endomorphism $h\colon A^{\Z}\to A^{\Z}$ takes it to itself. Furthermore, periodic points are dense in $A^\Z$. 
This allows to construct an inverse for $h$ defined on all of $A^{\Z}$.

\subsubsection*{Case 3: $G$ is residually finite}
For a normal subgroup $H\triangleleft G$ of finite index, define a $G/H$-periodic point of the shift $A^{G}$ as a map $x\colon G\to A$ that is constant on left $H$-cosets. (For instance, in case $G=\Z$ an $n$-periodic point is $\Z/n\Z$-periodic in this sense.) Now the argument from Case 2 applies.

\subsubsection*{Case 4: Another proof for $G=\Z$}$\,$
\\[1.5mm]
Let $X$ be a subshift of $A^{\Z}$, that is, a closed $G$-invariant subset. For every $n$ define $X\upharpoonright {A^{[-n,n]}}$ as the set of restrictions of all elements of $X$ to $[-n,n]$. The 
{\em topological entropy} of $X$ is defined by the formula
\begin{equation}
\label{eq:ht}
{\mathrm{ht}}(X)=\lim_{n\to\infty}
\frac{\log_{\abs A}\vert X\upharpoonright {A^{[-n,n]}}\vert}{n}.
\end{equation}
(Here one needs to work a little bit to show that the limit exists and equals the infimum, by proving and using the inequality $\abs{X_{n+m}}\leq \abs{X_n}\abs{X_m}$.) It is easy to see that ${\mathrm{ht}}(X)=1\iff X=A^\Z$: once the ratio on the r.h.s. goes below one, it cannot bounce back. A less obvious fact is that isomorphisms between subshifts of $A^G$ preserve the entropy value. Here one must use the observation that every such morphism admits a local representation in the sense that it can be fully recovered from a function defined on a suitable finite power of $A$. The combination of the two properties settles the case.
\\[3mm]
This generalizes to:
\subsubsection*{Case 5: $G$ is amenable} $\,$
\\[1.5mm]
The above observation about local representations of morphisms between subshifts can be used to show that surjunctivity of a group is a local property, and hence the argument carries over to 
\subsubsection*{Case 6: Initially subamenable groups}$\,$ \\[1.5mm]
The same idea can be stretched even further and made to work for 
\subsubsection*{Case 7: $G$ is sofic} $\,$
\\[1.5mm]
Here the definition of soficity is the original one by Gromov (that is, the property ({\large $\star$}) in Theorem \ref{th:graph}.)

For a highly readable presentation, Benjy Weiss's article \cite{weiss} is recommended.
 
\section{Von Neumann algebras and tracial ultraproducts}

Let $\H$ be a Hilbert space. Denote by $\B(\H)$ the $\ast$-algebra of all bounded operators on $\H$ equipped with the uniform norm. For example, if $\H=\C^n$, then $\B(\H)$ is the algebra $M_n(\C)$ of all $n\times n$ matrices with complex entries, equipped with the usual matrix addition, multiplication, and conjugate transpose. 

A $C^\ast${\em -algebra} is a Banach $\ast$-algebra isomorphic to a norm-closed subalgebra of $\B(\H)$. (See e.g. Weaver's survey \cite{weaver} for details.)

A {\em von Neumann algebra} is a unital $C^\ast$-algebra $M$ isomorphic to a weakly closed (equivalently: strongly closed) $C^\ast$-subalgebra of  $\B(\H)$. 
(The weak topology on $\B(\H)$ is the restriction of the Tychonoff product topology from the power $(\H_w)^{\H}$, where $\H_w$ is $\H$ equipped with its weak topology, while the strong topology is induced from ${\H}^{\H}$ where $\H$ carries the norm topology.) Equivalently, von Neumann algebras can be described as those $C^\ast$-algebras isometrically isomorphic to a dual space of some Banach space. This Banach space, called a predual of $M$ and denoted $M_\ast$, is necessarily unique up to isometric isomorphism. The $\sigma(M,M_\ast)$-topology on $M$ is called the {\em ultraweak topology,} and a $C^\ast$-algebra morphism between two von Neumann algebras is a von Neumann algebra morphism if it is ultraweak continuous.

\begin{example}
$\B(\H)$ itself is a von Neumann algebra.
\end{example}

\begin{example} If $G$ is a group, the {\em group von Neumann algebra} of $G$, denoted $VN(G)$, is the strong closure of a subalgebra of $\B(\ell^2(G))$ generated by all left translation operators $\lambda_g$, $g\in G$. 

Here $\lambda_g(f)(x)=f(g^{-1}x)$, while $\ell^2(G)$ denotes the Hilbert space of all $2$-summable complex-valued functions on $G$.
\end{example}

A von Neumann algebra $M$ is called a {\em factor} if the center of $M$ is trivial, that is, consists of only the constants, $\C\cdot 1$. 

\begin{example}
$\B(\H)$ is a factor. 
\end{example}

\begin{example}
If all the conjugacy classes of a group $G$ except $\{e\}$ are infinite (one says that $G$ has i.c.c. property, from ``infinite conjugacy classes''), then $VN(G)$ is a factor. For instance, every non-abelian free group, as well as the group $S_\infty^{fin}$ of all permutations of $\omega$ having finite support, has the i.c.c. property.
\end{example}
 
A von Neumann algebra $M$ is {\em approximately finite dimensional} (AFD) if it contains an increasing chain of finite dimensional subalgebras whose union is strongly dense in $M$.

\begin{example} Von Neumann algebras $VN(S_\infty^{fin})$ and $\B(\H)$ are AFD.
\end{example}

A {\em trace} on a $C^\ast$-algebra is a positive linear functional $\tr\colon A\to \C$ with $\tr(AB)=\tr(BA)$. 

\begin{example}
For every group $G$, the group von Neumann algebra $VN(G)$ has a trace, 
% If $G$ has i.c.c. property, then $VN(G)$ is a finite factor, where the trace
determined by the conditions $\tr(e)=1$ and $\tr(g)=0$, $g\neq e$.  
\end{example}

A von Neumann factor $M$ is {\em finite} if it admits a trace.

\begin{example}
$M_n(\C)$ is a finite factor, and the trace is the usual (normalized) trace of a matrix.
\end{example}

\begin{example}
If $G$ has the i.c.c. property, then $VN(G)$ is a finite factor.
\end{example}

A von Neumann factor is of type $II_1$ if it is finite and at the same time has infinite dimension as a Banach space. 

\begin{theorem}[Murray and von Neumann]
There exists a unique AFD von Neumann factor of type $II_1$ up to von Neumann algebra isomorphism.
\end{theorem}

This factor is denoted $R$ (not to be confused with the real line, $\R$).

\begin{example}
$R\cong VN(S_{\infty}^{fin})$.  
\end{example}

The normed space ultraproduct of a family of $C^\ast$-algebras is again a $C^\ast$-algebra in a natural way. This follows from submultiplicativity of the norm of a $C^\ast$-algebra, as well as from a characterization of $C^\ast$-algebras as those Banach algebras with involution satisfying the identity $\norm{x^\ast x}=\norm x^2$.

However, for von Neumann algebras an analogous statement is no longer true. Indeed, every von Neumann algebra is necessarily {\em monotonically complete} in the sense that every increasing bounded above net of positive elements has a least upper bound. (Cf. e.g. \cite{sakai}, Lemma 1.7.4.)
% (The $C^\ast$-algebras with this property are called {\em $AW^\ast$-algebras,} and they form a strictly intermediate class between $C^\ast$-algebras and von Neumann algebras.) 
But ultraproducts are known not to behave well with regard to order completeness. 

\begin{example}
\label{ex:counter}
Let $\mathcal U$ be a non-principal ultrafilter on the set of natural numbers. Denote by $[1,\nu]$ the ultraproduct modulo $\mathcal U$ of the family of totally ordered sets $[1,n]$, $n\in\N$, viewed as an ordered set. (Here we think of $\nu$ as an infinitely large integer.) 
The order structure of $[1,\nu]$ is well understood. If we denote by $\omega^\ast$ the set of positive integers with inverse order, then 
the order type of the segment $[1,\nu]$ is 
\[\omega\cup\eta\times  (\omega^\ast\cup\omega)\cup \omega^\ast,\]
where $\eta$ denotes a densely ordered set without the first and the last elements (which is, in fact, also countably saturated), and the order on the product is lexicographic. The leftmost copy of $\omega$ corresponds to the standard natural numbers sitting inside of the non-standard natural numbers as an initial segment. 

The $C^\ast$-algebra (normed space) ultraproduct of finite-dimensional commutative von Neumann algebras $\ell^{\infty}(n)$, $n\in\N$, with regard to $\mathcal U$ is easily seen to embed, in a canonical way, into the von Neumann algebra $\ell^\infty [1,\nu]$. 
For every $\xi\in [1,\nu]$ the characteristic function $\chi_{[1,\xi]}$ of the interval belongs to the ultraproduct. The sequence $(\chi_{[1,n]})_{n=1}^\infty$ is bounded, increasing, consists of positive elements, and yet has no least upper bound in $\left(\prod \ell^{\infty}(n)\right)_{\mathcal U}$. Indeed, the least upper bound of this sequence in the larger von Neumann algebra $\ell^\infty [1,\nu]$ is $\chi_{\omega}$, the characteristic function of the standard natural numbers, which is easily checked not to belong to the ultraproduct. (In the nonstandard analysis parlance, $\chi_{\omega}$ is an {\em external} function.) 
\end{example}

The above is more than an isolated counter-example: the same phenomenon will be observed in the $C^\ast$-algebra ultraproduct of any non-trivial family of von Neumann algebras.

Thus, for von Neumann algebras the ultraproduct construction needs to be modified. Example \ref{ex:counter} actually suggests how: we will need to factor out the members of the offending sequence $(\chi_{[1,n]})$. While each of them has $\ell^\infty$ norm one, their Hilbert-Schmidt norm, given by
\[\norm x_2 = \lim_{n\to{\mathcal U}}\tr_n(x^\ast x)^{1/2},\]
where $\tr_n(x)=\frac 1n\sum x_i$, vanishes as $n\to\infty$. This means that if we divide the $\ell^\infty$-direct sum of von Neumann algebras by a larger ideal of Hilbert-Schmidt norm infinitesimals, our counter-example will simply disappear inside of this ideal. At the same time, the Hilbert-Schmidt norm is not submultiplicative and is therefore unsuitable for forming the algebra of finite elements, for which we will still have to resort to the usual norm.

We will only present the construction in a particular case where all the $M_\alpha$ are factors of type $II_1$, equipped with traces normalized so that $\tr_\alpha(1)=1$. 

Introduce on every $M_\alpha$ the (normalized) Hilbert-Schmidt norm
\[\norm x_2 = \tr_\alpha(x^\ast x)^{1/2}.\]
Consider the $C^\ast$-algebra
\[\oplus^{\ell^\infty}_{\alpha\in A} M_\alpha/{\mathscr I},\]
where the $\ell^\infty$-direct sum is formed with regard to the {\em standard norms} on $M_\alpha$, while the infinitesimals are formed with regard to the {\em Hilbert-Schmidt norms:}
\[{\mathscr I}=\left\{x\colon \lim_{\alpha\to{\mathcal U}}\norm{x_\alpha}_2=0\right\}.\]
This $C^\ast$-algebra turns out to be a factor of type $II_1$, called the ({\em tracial}) {\em ultraproduct} of $M_\alpha$. The ideal $\mathscr I$ is larger than the ideal of norm infinitesimals, but it is not weakly closed in $\oplus^{\ell^\infty}_{\alpha\in A} M_\alpha$, so 
the result is quite surprising. 

However, the verification is not especially difficult. Here is an outline. Denote by $\H$ the Hilbert space completion of the inner product space $\oplus^{\ell^\infty}_{\alpha\in A} M_\alpha/{\mathscr I}$, where the inner product is defined by the trace $\tr(x)=\lim_{\alpha\to{\mathcal U}}\tr_\alpha(x_{\alpha})$. The algebra $M=\oplus^{\ell^\infty}_{\alpha\in A} M_\alpha/{\mathscr I}$ acts on itself by left multiplication, and this action extends by continuity to a faithful $C^\ast$-algebra representation of $M$ in $\H$. (This is the {\em GNS construction,} so named after Gelfand--Na\"\i mark and Segal.) 

The weak closure of $M$ in $\B(\H)$, which we will denote $\hat M$, is a von Neumann algebra. The trace extends by continuity over $\hat M$, and the corresponding Hilbert-Schmidt topology, while still Hausdorff, is coarser than the weak topology. The technique of polar decompositions (\cite{sakai}, Thm. 1.12.1) allows one to conclude that the unitary group of $M$ is the ultraproduct of the unitary groups $U(M_\alpha)$, equipped with their normalized Hilbert-Schmidt metrics, modulo $\mathcal U$. As we have mentioned before, the ultraproduct of a family of metric groups modulo a non-principal ultrafilter is a complete group, and thus closed whenever it embeds into a Hausdorff topological group as a topological subgroup. By \cite{pedersen}, 2.3.3, $U(M)$ is weakly dense in $U(\hat M)$, therefore Hilbert-Schmidt dense, and so the two groups coincide, and consequently $M=\hat M$ (because every element of $M$ is a linear combination of at most four unitaries). 

Let us finally show that $M$ is a factor. We will only do this in a particularly transparent case where $M_\alpha$ are matrix algebras $M_n(\C)$. 
Because of the last remark in the previous paragraph, it is enough to show that the center of the unitary group $U(M)$ is reduced to the circle group, because this will mean that the center of $M$ consists of scalars. In other words, we want to prove that every sequence $(u_n)$ of unitaries that is an approximate centralizer (that is, for each other sequence of unitaries $(v_n)$ one has $u_nv_nu_n^\ast v_n^\ast\overset{\mathcal U}\to 1$) converges to some $\lambda\cdot 1$, $\abs\lambda=1$. We prove this by contraposition. Let all the unitaries $(u_n)$ be at a Hilbert-Schmidt distance $\geq\e$ from the corresponding set of constants $\lambda\cdot 1_n$. 
In view of bi-invariance of the distance, this means, informally, that the eigenvalues of $u_n$ form a ``somewhat non-constant'' family. Let $w_n$ be a unitary that diagonalizes $u_n$. There is a coordinate permutation $\sigma_n$ with the property $\norm{w_nu_nw_n^\ast-\sigma_nw_nu_nw_n^\ast\sigma_n^{-1}}\geq \e$, implying that $(u_n)$ does not asymptotically centralize the sequence $(w_n^{\ast}\sigma_nw_n)$. 
% implies that the diameter of the orbit of each $u_n$ under conjugations, $vu_nv^\ast$, $v\in U(n)$, is $\geq\e$. 

The construction of tracial ultraproduct appears simultaneously in two articles independently published in 1970 by McDuff \cite{mcduff} and by Janssen \cite{janssen}.
A good presentation can also be found e.g. in Pisier's book \cite{pisier}, section 9.10. There is also a proof in \cite{BO}, Appendix A, which is however quite terse.

Notice that, assuming CH, all tracial ultraproducts of a fixed separable factor of type $II_1$ with regard to a nonprincipal ultrafilter on natural numbers are isomorphic between themselves. This was proved in \cite{GH}, essentially using the fact that ultraproducts of metric structures are countably saturated in an appropriate sense, as explained in \cite{BYBHU}. (Although this author must confess that Freiling's dart-throwing argument \cite{freiling} leaves him prejudiced against ever assuming the validity of the Continuum Hypothesis...)

\section{Connes' Embedding Conjecture\label{s:cec}}

Here is the most celebrated open problem of all those mentioned in this article. Recall that $R$ denotes the (unique) approximately finite dimensional factor of type $II_1$.

\begin{conjecture}[Connes' Embedding Conjecture] Every separable factor of type $II_1$ embeds into a suitable tracial ultrapower, $\left(R^{\N}\right)_{\mathcal U}$, of $R$.
\end{conjecture}

Connes himself proved the result for $VN(F_2)$ in \cite{connes-injective}.
% , that is, in our terminology, he has established that the free non-abelian group is hyperlinear. 
Independently and at the same time, this was also proved by Simon Wassermann (\cite{wassermann}, Lemma on p. 245). Connes then went on to remark (\cite{connes-injective}, p. 105): ``Apparently such an imbedding ought to exist for all $II_1$ factors...'' In the decades that followed the conjecture has become one of the central open problems of operator algebra theory. Through the work of Kirchberg (see e.g. \cite{kirchberg}), many equivalent forms of the conjecture have become known. For an in-depth discussion, see \cite{ozawa} and \cite{pisier}, section 9.10.

\begin{remark}
In operator algebra literature, the von Neumann algebra ultrapower $\left(R^{\N}\right)_{\mathcal U}$ is usually denoted $R^\omega$. 
The reader should beware of this notation. Here $R$ is {\em not} the real line $\R$ but the unique AFD factor of type $II_1$, and $\omega$ is {\em not} the first infinite ordinal number, but a (generic symbol for a) nonprincipal ultrafilter. The notation is traditional in operator algebra theory, and takes some getting used to for mathematicians coming from outside of the area.
\end{remark}

A particularly interesting case is that of group von Neumann factors.

\begin{conjecture}[Connes Embedding Conjecture for Groups]
\label{q:cecg?}
For every countable group $G$, the group von Neumann algebra $VN(G)$ embeds into a tracial ultrapower $R^\omega=\left(R^{\N}\right)_{\mathcal U}$ of $R$. 
\end{conjecture}

\begin{remark}
\label{r:R}
Hyperlinear groups can be alternatively characterized as subgroups of metric ultrapowers of $U(R)$, the unitary group of the AFD factor of type $II_1$ equipped with the Hilbert--Schmidt distance. This follows from Theorem \ref{th:onegroup}. For a relationship between $U(R)_2$ and the group $U(\infty)_2$ from Example \ref{ex:uinfty2}, as well as for more examples of related groups, see \cite{blattner} and also especially \cite{dlHP}.
% As it is well-known that the full group $[{\mathcal R}]$ of the tail equivalence relation is isomorphic to a topological subgroup of $U(R)$ with the strong operator topology (which happens to be the topology induced by the Hilbert-Schmidt distance), the two characterizations are in perfect agreement with each other.
\end{remark}

\begin{theorem}[Kirchberg \cite{kirchberg94}; Radulescu \cite{radulescu00}, Prop. 2.5; Ozawa \cite{ozawa}, Prop. 7.1]
Let $G$ be a countable group. Then $VN(G)$ embeds into $R^\omega$ if and only if $G$ is hyperlinear.
\label{th:iff}
\end{theorem}

\begin{proof}[Proof, sketch]
$\Rightarrow$: 
Suppose that $VN(G)$ embeds into $R^{\omega}$. As $G$ is contained in the unitary group $U(VN(G))$ as a subgroup, $G$ embeds into $U(R^{\omega})$ as a subgroup. As we noted elsewhere, the latter is isomorphic to a metric ultrapower of $U(R)$, and one concludes by Remark \ref{r:R}.
\\[4mm]
$\Leftarrow$: Let, as before, the tracial ultraproduct $R^{\omega}$ act on the Hilbert space $\H$ completion of $R^{\omega}$ equipped with the inner product $\langle x,y\rangle = \tr(x^{\ast}y)$, by assigning to every $x\in R^{\omega}$ the operator of left multiplication by $x$:
\[y\mapsto xy.\]
Now assume that a group $G$ embeds into the unitary group of $R^{\omega}$ as a subgroup. In view of Remark \ref{r:sqrt2}, 
one can assume that images of elements of $G$ are at the Hilbert-Schmidt distance $\sqrt 2$ from each other, that is, pairwise orthogonal in $\H$.
In other words, the restiction of the trace of $R^{\omega}$ to $G$ is Kronecker's delta $\delta_e$, and so $\H$ contains $\ell^2(G)$ as a Hilbert $G$-submodule. 
Denote by $M$ the weakly closed linear span of $G$ in $R^{\omega}$. The preceding sentence implies that there exists a canonical von Neumann algebra morphism from $M$ onto $VN(G)$. Since $R^{\omega}$ is a factor of type $II_1$, on the unitary group $U(R^\omega)$ the Hilbert-Schmidt topology determined by the trace coincides with the strong (and the weak) topologies given by the GNS representation (\cite{jones}, Proposition 9.1.1). Consequently, the same is true of $U(M)$, which implies that the restriction of the morphism $M\to VN(G)$ to $U(M)$ is in fact an isomorphism of topological groups (it is a Hilbert-Schmidt isometry). This leads to conclude that the morphism $M\to VN(G)$ has trivial kernel $N$, for otherwise the subgroup $(N+1)\cap M$ would be nontrivial as the unitary group of the unitalization of $N$.
\end{proof}

Thus, Connes' Conjecture for Groups (problem \ref{q:cecg?}) is equivalent to the statement that every group is hyperlinear (problem \ref{q:hyperlinear?}).

\begin{remark} For a countable group $G$ the following properties are equivalent.
\begin{enumerate}
\item $G$ embeds into the unitary group of the $C^\ast$-algebra ultraproduct of  matrix algebras $M_n(\C)$, $n\in\N$ with regard to some (any) nonprincipal ultrafilter on natural numbers.
\item $G$ embeds into the metric ultraproduct of unitary groups $U(n)$, $n\in\N$ formed with regard to the uniform operator metric.
\item $G$ embeds into a metric ultrapower of the group $U(\ell^2)_c$ of all compact perturbations of the identity, 
\[U(\ell^2)_c = \{u\in U(\ell^2)\colon u -{\mathbb{I}}\mbox{ is compact}\},\]
equipped with the uniform operator metric.
\end{enumerate}

The present author is unaware of any study of this class of groups, in particular, of the answers to the following questions.
\label{r:compact}
\end{remark}

\begin{question}
What is the relationship between the class of hyperlinear groups and the class of groups described in Rem. \ref{r:compact}? In particular, is every hyperlinear group contained in this class, and vice versa?
\end{question}

\begin{question}
Does the class of groups described in Rem. \ref{r:compact} contain every countable group?
\end{question}

\section{Some classes of groups to look at}

The two candidates for a counter-example are mentioned in Questions \ref{q:T?} and \ref{q:weiss}.

\begin{question}[Cf. Ozawa \cite{ozawa}]
Let $G$ be an infinite simple group with Kazhdan's property $(T)$. Can it be hyperlinear (sofic)?
\label{q:T?}
\end{question}

For theory of groups with property (T), we refer to \cite{dlHV} and especially \cite{BdlHV}. For a way to construct groups with a combination of properties mentioned in Problem \ref{q:T?}, see \cite{gromov87}.

Since in Gromov's construction the groups in question arise as direct limits of word hyperbolic groups, a positive answer to the following question would destroy any hope for a counter-example stemming from problem \ref{q:T?}.
A finitely-generated group $G$ with a set of generators $V$ is {\em word-hyperbolic} \cite{gromov87} if there is a constant $\delta>0$ with the property that for every three points $x,y,z$ the shortest path $[x,y]$ joining the vertices $x$ and $y$ is contained in the $\delta$-neighbourhood of the union of shortest paths $[y,z]\cup [z,x]$. (This property does not depend on the choice of a set of generators, only the value of the constant $\delta>0$ does.) For instance, the free groups are word-hyperbolic, while the free abelian groups on $m\geq 2$ generators are not. For a brief introduction, see \cite{BO}, section 5.3.

\begin{question}[Ozawa \cite{ozawa}]
Is every word-hyperbolic group \cite{gromov87} hyperlinear (sofic)? 
\label{q:hyp?}
\end{question}

This is, in fact, a weaker form of a famous open problem of combinatorial group theory:

\begin{question}
Is every word-hyperbolic group residually finite?
\label{q:res?}
\end{question}

\begin{question}[Weiss \cite{weiss}]
Is the free Burnside group of a finite exponent $n$ sofic? 
\label{q:weiss}
\end{question}

The free Burnside group of exponent $n$ is the quotient of the free group $F_m$ on $m\geq 2$ generators by the normal subgroup generated by the $n$-th powers of elements of $F_m$. As shown by Adian \cite{adian}, free Burnside groups of sufficiently large odd exponent are non-amenable. 

As the free Burnside group of large exponent is a direct limit of word hyperbolic groups (see \cite{olshanskii}, Ch. 6), a negative answer to Question \ref{q:weiss} would in particular imply the existence of a non residually finite word hyperbolic group, hence a negative answer to Question \ref{q:res?}. (This was pointed out to me by Goulnara Arzhantseva.)

A group $G$ has {\em Haagerup property} (or: is {\em a-$T$-menable}) (\cite{CCJJV}; also \cite{BO}, Section 12.2) if there is a sequence of positive definite functions $\phi_n$ on $G$, vanishing at infinity and converging pointwise to $1$. This property can be regarded as both a strong negation of Kazhdan's property $(T)$ and as a weak form of amenability: every amenable group has Haagerup property, but so does $F_2$. We suggest:

\begin{question}
Is every group with Haagerup property hyperlinear (sofic)? 
\end{question}

Here is another weak form of amenability. A countable discrete group $G$ is {\em amenable at infinity} (or {\em topologically amenable}) if for every finite subset $F\subseteq G$ and $\e>0$ there is a mapping $f$ from $G$ to the unit sphere of $\ell^2(G)$, having finite range and such that for all $g\in F$ and $x\in G$ one has $\norm{f(gx)-g\cdot f(x)}<\e$. (This equivalent description can be extracted from the results of \cite{HR}, and is nearly explicit in Proposition 4.4.5(2) from \cite{BO}.) For example, word-hyperbolic groups are amenable at infinity, see \cite{adams} and also Appendix B in \cite{ADR} (or else Theorem 5.3.15 in \cite{BO}).
It was shown that amenability at infinity is equivalent to a property of importance in operator algebra theory, exactness of a group (cf. \cite{HR,AD,ozawa01} or else section 5.1 in \cite{BO}).
The only known examples of (finitely generated) groups that are not amenable at infinity are those that can be obtained through Gromov's probabilistic method outlined in \cite{gromov03}.

\begin{question}
Is every group amenable at infinity ($=$ exact group)  hyperlinear (sofic)? 
\label{q:aminf?}
\end{question}

In view of \cite{guentner}, a positive answer to \ref{q:aminf?} would imply a positive answer to \ref{q:1rel?}.

\section{Equations in groups}

The following classical result can be put in direct connection with the present topic. Let $G$ be a group, $g_1,\ldots,g_n$ arbitrary elements of $G$, and $s_1,\ldots,s_n$ any integers such that $\sum s_i\neq 0$. Then the equation
\[x^{s_1}g_1x^{s_2}g_2\ldots x^{s_n}g_n=1\]
is called {\em regular}.

\begin{theorem}[Gerstenhaber and Rothaus \cite{GR}]
Every regular equation in a finite group $G$ has a solution in a finite group extending $G$. 
\label{th:gr}
\end{theorem}

The following is a long-standing conjecture that seems to be still open. (Cf. e.g. \cite{EJ} for some relatively recent advances.)

\begin{conjecture} Every regular equation in a group $G$ has a solution in some  group extending $G$.
\label{q:over?}
\end{conjecture}

There is an interesting link between this problem and the topic of the present survey. 
Gerstenhaber and Rothaus deduce their Theorem \ref{th:gr} from the following intermediate result which is of great interest on its own. 

\begin{theorem}[Gerstenhaber and Rothaus \cite{GR}, p. 1532]
Every regular equation with coefficients in $U(n)$ has a solution in $U(n)$.
\end{theorem}

(The theorem was in fact established for all compact connected Lie groups, and for certain systems of equations.)

One deduces easily:

\begin{corollary}
Every regular equation with coefficients in a hyperlinear group $G$ has a solution in a suitable hyperlinear group extending $G$ (which can be taken as the metric ultraproduct of unitary groups of finite rank).
\end{corollary}

Thus, Connes' Embedding Conjecture for groups, if proved, would imply conjecture \ref{q:over?} on regular equations in groups, while any counter-example to the latter will disprove Connes' Embedding Conjecture. 

\section{Varia}

The present short survey, or rather a collection of introductory remarks, is emphatically not exhaustive. Here are some pointers to what has been left out.

A number of known conjectures have been settled in the positive for sofic groups, and so proving that every group is sofic would settle those conjectures as well. These include Kaplansky's Direct Finiteness conjecture \cite{ES0}, the Determinant Conjecture \cite{ES}, and some others \cite{thom}.

Sofic groups have been linked to theory of stochastic processes in infinite networks \cite{AL} and to cellular automation \cite{CSC}, and have been shown to admit a classification of their Bernoulli actions \cite{bowen}.

The known permanence properties of sofic groups are discussed in \cite{ES2}, see also \cite{GG}. For a discussion and references to similar properties of hyperlinear groups, see \cite{ozawa}.

Pierre de la Harpe has brought to my attention the following. Every (locally compact) group has a unique maximal (closed) amenable normal subgroup, called the amenable radical of $G$. (Day  \cite{day}, see also \cite{zimmer-book}, Proposition 4.1.12.) 

\begin{question}[Pierre de la Harpe]
Can one define the sofic radical or the hyperlinear radical of a group?
\end{question}

This requires proving the following permanence property: if $N_1,N_2\triangleleft G$ are two normal sofic (respecively hyperlinear) subgroups of $G$, then the group $N_1N_2$ is sofic (resp., hyperlinear). 
% The concept will become of interest after (and if) the existence of non-sofic and/or non-hyperlinear groups is established. 

% Every countable discrete group $G$ has
% There is a well-established notion of an amenable radical due to Day (see e.g. \cite{zimmer-book}, Proposition 4.1.12). The amenable radical of a group is the inter
% As noted by Pierre de la Harpe (personal communication), one could define the sofic radical

Glebsky and Rivera Mart\'\i nez \cite{GM} introduce {\em weakly sofic groups,} as subgroups of metric ultraproducts of arbitrary families of finite groups equipped with bi-invariant metrics. They put forward the following conjecture, linking it to the deep investigation of Herwig and Lascar \cite{HL}:

\begin{conjecture}[\cite{GM}] Every group is weakly sofic.
\end{conjecture}

And finally, the following comment was made by Ilijas Farah when explaining to me the contents of his talk \cite{farah}. 

\begin{remark}
Set theory suggests that the solutions to the main problems of this survey (such as questions \ref{q:sofic?} and \ref{q:hyperlinear?}) are unlikely to be independent -- at the very least, they will not be affected by adding to the axioms of ZFC some of the most popular additional axioms, such as the Continuum Hypothesis (CH), Martin's Axiom (MA), or the Axiom of Constructibility ($V=L$).
The reason is that the statements of corresponding conjectures (denote them $\phi$) are absolute between transitive models of ZF containing all countable ordinals, and consequently the following are equivalent:
\begin{enumerate}
\item ZFC $\vdash \phi$
\item ZFC+CH $\vdash \phi$
\item ZFC+MA $\vdash \phi$
\item ZFC+ (V=L) $\vdash \phi$
\item ZFC+A $\vdash \phi$
\end{enumerate}
Here $A$ denotes any axiom of ZFC such that the constructible unverse, L, has
a forcing extension in which $A$ holds. Moreover, `ZFC' can
be replaced with `ZF' (which is most interesting in (1)). Cf. Theorem 13.15 on p. 175, and remark following it, in \cite{kanamori}.
\end{remark}

\section{Some reading suggestions}

A very good introduction to the subject of sofic groups is Weiss' survey \cite{weiss}, which should be followed by the paper \cite{ES} of Elek and Sz\'abo (treating both sofic and hyperlinear groups). Section 9.10 of Pisier's book \cite{pisier}, devoted to ultraproducts and Connes' Conjecture, makes for an enjoyable (and more or less self-contained) read. It can be complemented by a recent introduction by Weaver \cite{weaver} to the theory of $C^\ast$-algebras and von Neumann algebras for logicians. 
Ozawa's survey \cite{ozawa} paints a broad picture of Connes' Embedding Conjecture in its many equivalent forms and ramifications, and the newly-published book \cite{BO} of Brown and Ozawa is an even more comprehensive source. After that, one cannot ignore seminal articles by Connes \cite{connes-injective} and Gromov \cite{gromov99} --- difficult but stimulating reads.

\section*{Conclusion}

Answers to some of the questions discussed in this article seem to be suggested by the following dichotomy attributed to Gromov (cf. \cite{ozawa}):

\begin{quote}
{\em
Any statement about all countable groups is either trivial or false.}
\end{quote}

In particular, in the spirit of this dictum, one needs to look for counter-examples to Connes' Embeddings Conjecture for Groups --- unless a proof of the conjecture turns out to be unexpectedly simple.

\section*{Acknowledgements}
These notes are based on lectures given at the CIRM Luminy workshop {\em Petit groupe de travail sur la conjecture de Connes} (5--7 Dec. 2005) and the Fields Institute Set Theory Seminar (Jan. 2007).
I am thankful to Emmanuel Germain and to Stevo Todorcevic respectively for invitations to speak at the above forums. My gratitude extends  to Matthew Mazowita who brought my attention to the paper \cite{ES}, as well as to Sebastian Dewhurst and other participants of a summer 2006 student research seminar in Ottawa where we went through many proofs in detail.
At some point, Ilijas Farah had convinced me that writing my lecture up might be of interest to logicians. Ilijas has suggested a large number of improvements to every single version of the article, with further suggestions coming from Goulnara Arzhantseva, Lewis Bowen, Beno\^\i t Collins, Pierre de la Harpe, Gab\'or Elek, Thierry Giordano, Alexander Kechris, Denis Osin, Simon Thomas, two anonymous referees, and an anonymous senior logician. I am grateful to all of them. A conductive environment for discussions was provided by
the Fields Institute workshop around Connes' Embedding Problem held from May 16--18, 2008 at the University of Ottawa. My work was supported by NSERC operating grants (2003-07, 07--) and University of Ottawa internal grants.


\begin{thebibliography}{100}

\bibitem{adams} S. Adams, {\em Boundary amenability for word hyperbolic groups and an application to smooth 
dynamics of simple groups,} Topology \textbf{33} (1994), 765-783. 

\bibitem{adian}  S. I. Adian, {\em Random walks on free periodic groups,} Izv. Akad. Nauk SSSR, Ser. Mat. \textbf{46} (1982), 1139--1149.

\bibitem{AL} D. Aldous and R. Lyons, {\em Processes on unimodular random networks,} Electron. J. Probab. \textbf{12} (2007), 1454--1508.

\bibitem{AD} C. Anantharaman-Delaroche, {\em Amenability and exactness for dynamical systems and their $C\sp *$-algebras,} Trans. Amer. Math. Soc. \textbf{354} (2002), 4153--4178.

\bibitem{ADR} C. Anantharaman-Delaroche and J. Renault, 
{\em Amenable groupoids.}
With a foreword by Georges Skandalis and Appendix B by E. Germain. Monographies de L'Enseignement Math\'ematique  \textbf{36}, L'Enseignement Math\'ematique, Geneva, 2000.

\bibitem{Bek1} M.E.B. Bekka,
{\it Amenable unitary representations of locally compact groups,}
Invent. Math. {\bf 100} (1990), 383--401.

\bibitem{BdlHV} M.B. Bekka, P. de la Harpe, and A. Valette,
{\em Kazhdan's Property $(T)$,} New Mathematical Monographs \textbf{11}, Cambridge University Press, 2008.

\bibitem{BS} J.L. Bell and A.B. Slomson, {\em Models and Ultraproducts: An Introduction,} North-Holland, Amsterdam--London, 2nd revised printing, 1971.

\bibitem{BYBHU} I. Ben Yaacov, A. Berenstein, C.W. Henson, and A. Usvyatsov, {\em Model Theory for Metric Structures,} 112 pages, to be published in a Newton Institute volume in the Lecture Notes series of the London Math. Society, current version available from {\tt http://www.math.uiuc.edu/\%7Ehenson/cfo/mtfms.pdf}

\bibitem{bergman}
G.M. Bergman, {\em Generating infinite symmetric groups,}
Bull. London Math. Soc. \textbf{38} (2006), 429--440. 

\bibitem{blattner} 
R.J. Blattner, {\em Automorphic group representations,} Pacific J. Math. \textbf{8} (1958), 665--677.

\bibitem{bowen} L. Bowen, {\em Isomorphism invariants for actions of sofic groups,} arXiv:0804.3582v1 [math.DS], 22 pages.

\bibitem{BO} N.P. Brown and N. Ozawa, {\em ${C}^*$-Algebras and Finite-Dimensional Approximations,} Graduate Studies in Mathematics \textbf{88},  American Mathematical Society, Providence, R.I., 2008.

\bibitem{BM} M. Burger and S. Mozes, 
{\em Finitely presented simple groups and products of trees,}
C. R. Acad. Sci. Paris S\'er. I Math. \textbf{324} (1997), 747--752. 

\bibitem{CSC} T. Ceccherini-Silberstein and M. Coornaert, 
{\em Injective linear cellular automata and sofic groups,}
Israel J. Math. \textbf{161} (2007), 1--15. 

\bibitem{CCJJV}
P.-A. Cherix, M. Cowling, P. Jolissaint, P. Julg, and A. Valette, 
{\em Groups with the Haagerup property.
Gromov's a-$T$-menability,} Progress in Mathematics, \textbf{197}, Birkh\"auser Verlag, Basel, 2001.

\bibitem{connes-injective} A. Connes, 
{\it Classification of injective factors,}
Ann. of Math. {\bf 104} (1976), 73--115.

\bibitem{DCK} D. Dacunha-Castelle and J.-L. Krivine, {\em Applications des ultraproduits \`a l'\'etude des espaces et des alg\`ebres 
de Banach,} Studia Mathematicae \textbf{41} (1972), pp. 315-334. 

\bibitem{day}
M. M. Day, {\em Amenable semigroups}, Illinois J. Math. \textbf{1} (1957), 509-Ð544. 

\bibitem{dlH3} P. de la Harpe, \textit{Topics in Geometric Group Theory,}
Chicago Lectures in Mathematics, University of Chicago Press, Chicago, IL, 2000. 

\bibitem{dlHP} P. de la Harpe and R.J. Plymen, {\em Automorphic group representations: a new proof of Blattner's theorem,} J. London Math. Soc. (2) \textbf{19} (1979), 509--522.

\bibitem{dlHV} P. de la Harpe and A. Valette,
{\it La propri\'et\'e {\rm (T)} de Kazhdan pour les groupes
localements compacts,} Ast\'erisque {\bf 175} (1989), 158 pp.

\bibitem{Dye} H.A. Dye, \textit{On groups of measure preserving transformations,
I, II}, Amer. J. Math. \textbf{81} (1959), 119--159; \textbf{85} (1963),
551--576.

\bibitem{EJ}
M. Edjvetand and A. Juh\'asz, {\em On equations over groups,}
Internat. J. Algebra Comput. \textbf{4} (1994), 451-468. 

\bibitem{ES0} G. Elek and E. Szab\'o, {\em Sofic groups and direct finiteness,}  J. Algebra \textbf{280} (2004), 426--434. 

\bibitem{ES}
G. Elek and E. Szab\'o, {\em
Hyperlinearity, essentially free actions and $L\sp 2$-invariants. The sofic property,}
Math. Ann. 332 (2005), no. 2, 421--441. 

\bibitem{ES2} G. Elek and E. Szab\'o, {\em On sofic groups,} J. Group Theory \textbf{9} (2006), 161--171. 

\bibitem{farah1} I. Farah, {\em The commutant of $B(H)$ in its ultrapower,} abstract of talk at the Second Canada--France Congress (1--5 June 2008, UQAM, Montr\'eal),\\ {\tt
http://www.cms.math.ca/Events/summer08/abs/set.html\#if}

\bibitem{farah} I. Farah, {\em Incompleteness, independence, and absoluteness or: When to call a set theorist?}, invited lecture at the Fields Institute Workshop Around Connes' Embedding Problem, University of Ottawa, May 16--18, 2008.\\
{\tt http://www.fields.utoronto.ca/programs/scientific/07-08/embedding/abstracts.html}

% \bibitem{FPS}
% I. Farah, N.C. Phillips, and J. Stepr\=ans,
% {\em The commutant of $L(H)$ in its ultrapower may or may not be trivial,}
% preprint, June 2008, 15 pp.

\bibitem{freiling} C. Freiling, {\em
Axioms of symmetry: throwing darts at the real number line,}
J. Symbolic Logic \textbf{51} (1986), 190--200. 

\bibitem{GH} 
L. Ge and D. Hadwin, {\em Ultraproducts of $C\sp *$-algebras,} in: Recent advances in operator theory and related topics (Szeged, 1999), 
Oper. Theory Adv. Appl., \textbf{127}, Birkh\"auser, Basel, 2001, 305-326.

\bibitem{GR}
M. Gerstenhaber and O.S. Rothaus, 
{\em The solution of sets of equations in groups,}
Proc. Nat. Acad. Sci. U.S.A. \textbf{48} (1962), 1531--1533. 

\bibitem{GG} L.Yu. Glebsky and E.I. Gordon,
{\em On surjunctivity of the transition functions of cellular automata on groups,} Taiwanese J. Math. \textbf{9} (2005), 511--520.

\bibitem{GM} L. Glebsky and L.M. Rivera Mart\'\i nez, {\em Sofic groups and profinite topology on free groups,} arXiv:math0709.0026, Feb. 2008, 5 pp.

\bibitem{gottschalk} 
W. Gottschalk, {\em Some general dynamical notions,} in: Recent Advances in Topological Dynamics, Lecture Notes Math. \textbf{318}, Springer-Verlag, Berlin, 1973, pp. 120--125.

\bibitem{Gre1} F.P. Greenleaf,
{\it Invariant Means on Topological Groups,}
Van Nostrand Mathematical Studies {\bf 16},
Van Nostrand -- Reinhold Co., NY--Toronto--London--Melbourne, 1969.

\bibitem{gromov87} M. Gromov, {\em Hyperbolic groups,} in: Gersten, S.M. (ed.) Essays in group theory (M.S.R.I. Publ. \textbf{8}, Springer 1987, pp. 75--263.

\bibitem{gromov99}
M. Gromov, {\em Endomorphisms of symbolic algebraic varieties,} J. Eur. Math. Soc. (JEMS) 1 (1999), no. 2, 109--197.

\bibitem{gromov03} M. Gromov, {\em Random walk in random groups,} Geom. Funct. Anal. \textbf{13} (2003), 73--146.

\bibitem{guentner}
E. Guentner, {\em Exactness of the one relator groups,}
Proc. Amer. Math. Soc. \textbf{130} (2002), 1087--1093.

\bibitem{HI} C.W. Henson and J. Iovino,
\textit{Ultraproducts in analysis}, in: {Analysis and Logic (Mons, 1997)},
{London Math. Soc. Lecture Note Ser.}, \textbf{262}, {1--110},
{Cambridge Univ. Press}, {Cambridge}, {2002}.

\bibitem{HL}
B. Herwig and D. Lascar, {\em
Extending partial automorphisms and the profinite topology on free groups,}
Trans. Amer. Math. Soc. \textbf{352} (2000), 1985--2021. 

\bibitem{hewitt}
E. Hewitt, {\em Rings of real-valued continuous functions. I,}
Trans. Amer. Math. Soc. \textbf{64} (1948), 45--99. 

\bibitem{HR}
N. Higson and J. Roe, 
{\em Amenable group actions and the Novikov conjecture,}
J. Reine Angew. Math. \textbf{519} (2000), 143--153. 

\bibitem{hjorth}
G. Hjorth, {\em An oscillation theorem for groups of isometries,} Geometric and Functional Analysis \textbf{18} (2008), 489--521.

\bibitem{HM}
G. Hjorth and M. Molberg, {\em
Free continuous actions on zero-dimensional spaces,}
Topology Appl. \textbf{153} (2006), 1116--1131. 

\bibitem{janssen}
G. Janssen, {\em Restricted ultraproducts of finite von Neumann algebras,} in:  Contributions to Non-Standard Analysis (Sympos., Oberwolfach, 1970), Studies in Logic and Found. Math., Vol. 69, North-Holland, Amsterdam, 1972, pp. 101--114. 

\bibitem{jones} V.F.R. Jones, {\em Von Neumann Algebras,} lecture notes, version of 13th May 2003, 121 pp. Available from {\tt available from http://www.imsc.res.in/$\sim$vpgupta/hnotes/notes-jones.pdf}

\bibitem{kanamori}
A. Kanamori, {\em The Higher Infinite: Large Cardinals in Set Theory from Their
Beginnings}, Springer, Berlin--Heidelberg--New York, Perspectives in Mathematical Logic, 1995.

\bibitem{KM} A.S. Kechris and B.D. Miller, \textit{Topics in Orbit Equivalence,}
Lect. Notes Math. \textbf{1852}, Springer--Verlag, Berlin, Heidelberg, 2004.

\bibitem{KR} A.S. Kechris and C. Rosendal, {\em Turbulence, amalgamation, and generic automorphisms of homogeneous structures,} Proc. Lond. Math. Soc. (3) \textbf{94} (2007), 302--350. 

\bibitem{kirchberg}  E. Kirchberg, \textit{On nonsemisplit extensions, 
tensor products and exactness of group $C\sp *$-algebras,}
Invent. Math.  \textbf{112} (1993), 449--489. 

\bibitem{kirchberg94} E. Kirchberg, {\em Discrete groups with Kazhdan's property ${\rm T}$ and factorization property are residually finite,} Math. Ann. \textbf{299} (1994), 551--563. 

\bibitem{kirchberg06} E. Kirchberg, {\em Central sequences in $C\sp *$-algebras and strongly purely infinite algebras,} in: Operator Algebras: The Abel Symposium 2004, Abel Symp., \textbf{1}, Springer, Berlin, 2006, pp. 175--231.

\bibitem{kropholler}
P. Kropholler, {\em Baumslag-Solitar groups and some other groups of cohomological dimension two,}
Comment. Math. Helvetici \textbf{65} (1990), 547--558.

\bibitem{los} J. \L o\'s, {\em Quelques remarques, th\'eor\`emes et probl\`emes sur les classes d\'efinissables d'alg\`ebres,} in: Mathematical interpretation of formal systems, North-Holland Publishing Co., Amsterdam, 1955,  pp. 98--113.

\bibitem{luxemburg} W.A.J. Luxemburg, {\em A general theory of monads,} in Applications of Model Theory to Algebra, Analysis, 
and Probability, Holt, Rinehart and Winston (New York, 1969), pp. 18--86.

\bibitem{MKS}  W. Magnus, A. Karrass, and D. Solitar, 
\textit{Combinatorial Group Theory,} reprint of the 1976 second edition,
Dover Publications, Mineola, NY, 2004.

\bibitem{malcev} A.I. Malcev, \textit{On isomorphic matrix representations
of infinite groups,} Math. Sbornik (N.S.) \textbf{8} (1940), no. 50, 
405--422.

\bibitem{mcduff} D. McDuff, {\em Central sequences and the hyperfinite factor,} Proc. London Math. Soc. (3) \textbf{21} (1970), 443--461.

\bibitem{NVT} L. Nguyen Van Th\'e, \textit{Structural Ramsey Theory of Metric Spaces and Topological Dynamics of Isometry Groups}, to appear as a volume of Memoirs of the American Mathematical Society.
%, based on a Ph.D. thesis available at {\tt http://tel.archives-ouvertes.fr/tel-00139239/fr/}

\bibitem{olshanskii} A.Yu. Ol'shanski\u\i, 
{\em Geometry of Defining Relations in Groups,} Mathematics and its Applications (Soviet Series), \textbf{70}, Kluwer, Dordrecht, 1991. 

\bibitem{ozawa01} N. Ozawa, 
{\em Amenable actions and exactness for discrete groups,}
C. R. Acad. Sci. Paris S\'er. I Math. \textbf{330} (2000), 691--695. 

\bibitem{ozawa}  N. Ozawa, \textit{About the QWEP conjecture,}
Internat. J. Math.  \textbf{15} (2004), 501--530. 

\bibitem{Pat88} A.T. Paterson,
{\it Amenability,} Math. Surveys and Monographs
{\bf 29}, Amer. Math. Soc., Providence, RI, 1988.

\bibitem{pedersen}  G.K. Pedersen, {\em $C\sp{*} $-algebras and Their Automorphism Groups,} London Mathematical Society Monographs, \textbf{14}, Academic Press, Inc., London-New York, 1979. 

\bibitem{pisier} G. Pisier, \textit{Introduction to Operator Space Theory,}
London Math. Soc. Lect. Note Series, \textbf{294},
Cambridge Univ. Press, Cambridge, 2003. 

\bibitem{radulescu00}
F. Radulescu, {\em The von Neumann algebra of the non-residually finite Baumslag group $\langle a,b | a b^3 a^{-1} = b^2\rangle$ embeds into $R^{\omega}$}, arXiv:math/0004172v3, 2000, 16 pp.

\bibitem{robinson} A. Robinson, {\em Non-standard Analysis,} Studies in Logic and Foundations of Mathematics, North-Holland, Amsterdam--London, 1966 (second printing, 1970).

\bibitem{RS} C. Rosendal and S. Solecki, 
{\em Automatic continuity of homomorphisms and fixed points on metric compacta,}
Israel J. Math. \textbf{162} (2007), 349--371. 

\bibitem{sakai} S. Sakai, {\it $C^\ast$-Algebras and $W^\ast$-Algebras,}
Springer-Verlag, Berlin--Neidelberg--NY, 1971; Reprinted,
Springer, 1998.

\bibitem{sanov} I.N. Sanov, {\em A property of a representation of a free group,} Dokl. Akad. Nauk SSSR \textbf{57} (1947), 657--659 (in Russian).

\bibitem{thom} A. Thom, {\em Sofic groups and diophantine approximation,}
arXiv:math/0701294v3 [math.FA], 15 pages, Jan. 2007.

\bibitem{valette}
A. Valette, {\em Amenable representations and finite injective von Neumann algebras,} Proc. Amer. Math. Soc. \textbf{125} (1997), 1841--1843. 

\bibitem{VG} A.M. Vershik and E.I. Gordon, {\em Groups that are locally embeddable in the class of finite groups,} St. Petersburg Math. J. \textbf{9} (1998), 49--67.

\bibitem{wagon} S. Wagon, {\it The Banach--Tarski Paradox,}
Cambridge University Press, 1985.

\bibitem{wassermann} S. Wassermann, {\em On tensor products of certain group $C\sp{*} $-algebras,} J. Functional Analysis \textbf{23} (1976), 239--254. 

\bibitem{weaver} N. Weaver, {\em Set theory and $C^\ast$-algebras,}
The Bulletin of Symbolic Logic \textbf{13} (2007), 1--20.

\bibitem{weiss}
B. Weiss, 
{\em Sofic groups and dynamical systems,}
Sankhy\=a Ser. A 62 (2000), no. 3, 350--359. \\
Available at: {\tt
http://202.54.54.147/search/62a3/eh06fnl.pdf}

\bibitem{zimmer-book}
R.J. Zimmer, {\em Ergodic Theory and Semisimple Groups,} Monographs in Mathematics, \textbf{81}, BirkhŠuser Verlag, Basel, 1984.

\end{thebibliography}
\end{document}